\documentclass[leqno,11pt]{amsart}
\usepackage{amsfonts,latexsym}
\usepackage{amsmath}
\usepackage{amscd}
\usepackage{amsmath,amssymb,bm,graphics}
\usepackage[dvips]{graphicx}

\newcommand{\nequation}{\setcounter{equation}{0}}

\newcommand{\R}{{\mathbb R}}

\newcommand{\C}{{\mathbb C}}
\newcommand{\Z}{{\mathbb Z}}

\newtheorem{theorem}{Theorem}[section]
\newtheorem{proposition}[theorem]{Proposition}
\newtheorem{lemma}[theorem]{Lemma}
\newtheorem{corollary}[theorem]{Corollary}

\newtheorem{remark}[theorem]{Remark}

\newcommand{\proofbegin}{\noindent{\it Proof\,\,}}
\newcommand{\proofend}{\hfill$\Box$\bigskip}
\newenvironment{ackn}{\medskip \noindent \small
{\sl Acknowledgments.}}{\bigskip}

\newcommand{\marginnote}[1]
{
}
\newcounter{jl}

\newcounter{gm}

\newcounter{bk}

\input epsf
\date{April 7, 2008}
\title[Generalized Hunter-Saxton equation and the 
diffeomorphism group]
{Generalized Hunter-Saxton equation and \\the geometry of the group 
of circle diffeomorphisms}
\author{Boris Khesin}
\address{B.K.: Department of Mathematics, University of Toronto, ON M5S 2E4, Canada} 
\email{khesin@math.toronto.edu} 

\author{Jonatan Lenells}
\address{J.L.: Department of Applied Mathematics and Theoretical Physics, 
University of Cambridge, Cambridge CB3 0WA, UK}
\email{J.Lenells@damtp.cam.ac.uk} 

\author{Gerard Misio\l ek}
\address{G.M.: Department of Mathematics, University of Notre Dame, 
IN 46556, USA}
\email{gmisiole@nd.edu} 

\begin{document}
\maketitle

\begin{abstract}
We study an equation lying `mid-way' between the periodic Hunter-Saxton 
and Camassa-Holm equations, and which describes 
evolution of rotators in liquid crystals with external magnetic field and
self-interaction.
We prove that it is an Euler equation on the diffeomorphism group 
of the circle corresponding to a natural right-invariant Sobolev metric. 
We show that the equation is bihamiltonian and admits both cusped 
as well as smooth traveling-wave solutions which are natural candidates 
for solitons. 
We also prove that it is locally well-posed and establish results on 
the lifespan of its solutions. 
Throughout the paper we argue that despite similarities to 
the KdV, CH and HS equations, 
the new equation manifests several distinctive features 
that set it apart from the other three. 
\end{abstract}
\tableofcontents


\section{Introduction}\nequation

The equation we study in this paper is 
\begin{equation}\label{muHS}\tag{$\mu$HS}
- u_{txx} = - 2\mu(u) u_x + 2 u_x u_{xx} + u u_{xxx} 
\end{equation}
where $u=u(t,x)$ is a time-dependent function on the unit circle 
$S^1 = \mathbb{R}/\mathbb{Z}$ 
and $\mu(u) = \int_{S^1} u \, dx$ denotes its mean. 
The closest relatives of this equation are the Camassa-Holm equation
\begin{equation}\label{CH}\tag{CH}
u_t-u_{txx}+3uu_x = 2u_xu_{xx}+uu_{xxx},
\end{equation}
and the Hunter-Saxton equation
\begin{equation*}\label{HS}\tag{HS}
-u_{txx} = 2u_xu_{xx} + uu_{xxx},
\end{equation*}
both of which have recently attracted a lot of attention 
among the integrable systems and the PDE communities.
The Camassa-Holm equation, originally derived in \cite{F-F} as an abstract 
equation with a bihamiltonian structure, was introduced independently 
in \cite{C-H} as a shallow water approximation\footnote{For an interesting 
discussion of the attendant physical principles see \cite{J1}.}. 
The Hunter-Saxton equation first appeared in \cite{H-S} as an asymptotic 
equation for rotators in liquid crystals.
We show that (\ref{muHS}) can be viewed as a natural generalization
of the rotator equation once we allow interactions of rotators
and an external magnetic field. 

In order to compare equation (\ref{muHS}) with the other two equations 
we first recall a few facts about the latter. It is known that for 
spatially periodic functions (\ref{CH}) 
is the geodesic equation on the infinite-dimensional group 
$\mathcal{D}^s(S^1)$ of orientation-preserving diffeomorphisms 
of the unit circle $S^1 \simeq \R/\Z$ of (sufficiently high) 
Sobolev class $H^s$ and endowed with a right-invariant metric given 
at the identity by the $H^1$ inner product, see \cite{Mis98}, 
$$
\langle u, v \rangle_{H_1} = \int_{S^1} \bigl(uv + u_x v_x\bigr) dx.
$$ 
Equation (\ref{CH}) is a completely integrable system with a bihamiltonian 
structure 
and hence it possesses an infinite family of commuting Hamiltonian flows, 
as well as an associated sequence of conservation laws, 
see e.g. \cite{C-M} for the periodic case. 
Furthermore, it admits soliton-like solutions (called peakons) in both 
the periodic and the non-periodic setting, see \cite{C-H}. 

Equation (\ref{HS}), on the other hand, describes the geodesic flow 
on the homogeneous space of the group $\mathcal{D}^s(S^1)$ 
modulo the subgroup of rigid rotations $Rot(S^1) \simeq S^1$
equipped with the $\dot{H}^1$ right-invariant metric, which at 
the identity is
\begin{equation*}
\langle u, v \rangle_{\dot{H}_1} = \int_{S^1} u_x v_x dx.
\end{equation*}
This equation possesses a bihamiltonian structure as well, 
see e.g. \cite{H-Z, K-M}.

In Section \ref{GeoSec} we will show that the $\mu$HS equation 
has a similar geometric origin:
it describes the geodesic flow on $\mathcal{D}^s(S^1)$ with 
the right-invariant metric given at the identity by the inner product 
\begin{equation*} \label{mu-met} 
\langle u, v \rangle = \mu(u)\mu(v) + \int_{S^1} u_x v_x dx.
\end{equation*} 
(see also Section \ref{WP}). 
One motivation for considering this inner product 
comes from the observation made in \cite{L} that there is 
a natural isometry between $\mathcal{D}^s(S^1)/S^1$ equipped with 
the $\dot{H}^1$ metric 
and an open subset $U$ of the (infinite-dimensional) unit sphere 
in $L^2(S^1)$. 
Consequently, the group $\mathcal{D}^s(S^1)$ can be viewed as 
an $S^1$-bundle over $U$ with 
$\langle \cdot, \cdot \rangle$ 
inducing the simplest possible metric in this bundle which 
projects to the $\dot{H}^1$ metric on the base $U$. 
(In a sense, this picture can be thought of as an infinite-dimensional
version of the Hopf $S^1$-fibration.) 
Thus equation (\ref{muHS}) as a geodesic equation for 
$\mathcal{D}^s(S^1)$ 
with this metric is expected to have the `nicest' properties 
among the integrable equations related to this group,
owing to constant sectional curvature of the base $U$ and 
the `constant' metric along the fibers $S^1$.

In Section \ref{varphyssec} we will present two variational principles 
for the $\mu$HS equation and show how it arises as an asymptotic rotator 
equation in a liquid crystal with a preferred direction 
if we take into account the reciprocal action of dipoles on themselves.
The bihamiltonian structure of (\ref{muHS}) together with its attendant 
infinite hierarchy of conservation laws will be described in 
Sections \ref{Biham} and \ref{Virasoro}. 

In spite of the presence of the nonlocal expression $\mu(u)$ 
the $\mu$HS equation is a well-defined PDE (rather than a mixed 
integro-differential equation)\footnote{We mention here that a related 
(nonperiodic) equation obtained by replacing $\mu(u)$ by 
a constant has appeared previously as 
a special case of a model derived in \cite{M-N} for 
the propagation of surface waves in an ideal fluid with surface tension, 
and whose soliton solutions were studied using tools of algebraic 
and complex geometry in \cite{ACHM}.}. 
This is because the mean of any solution $u$ is constant in time 
and hence is completely determined by the initial condition.
In light of this fact the term $2\mu(u) u_x$ may appear rather harmless 
and might be expected to have only a minor influence on the phenomenology 
of the equation. 
However, in Section \ref{WP} we shall see that it has a strong effect 
on well-posedness of the associated Cauchy problem in that it is 
responsible for (\ref{muHS}) admitting global (in time) solutions. 
In the same section a blow-up mechanism for $\mu$HS solutions 
will also be described. 
In this respect the $\mu$HS equation resembles (\ref{CH}) 
rather than (\ref{HS}) as all classical solutions of the periodic HS equation 
except space-independent solutions blow up in finite time, 
see \cite{g} or \cite{Yin}. 
Let us point out that just like the breaking waves of (\ref{CH}) can be 
continued as global weak solutions (see \cite{B-C2}) the solutions 
of (\ref{HS}) can also be extended to all times (see \cite{B-C} 
or \cite{L19} for the periodic case).

Furthermore, whereas (\ref{HS}) does not have any bounded 
traveling-wave solutions at all, the $\mu$HS equation will be shown 
in Section \ref{waves} to admit traveling waves (in the periodic case) 
that can be regarded as the appropriate candidates for 
solitons.\footnote{That (\ref{HS}) does not admit 
any bounded traveling wave 
solutions is  related to the fact that the equation itself 
is defined on the homogeneous space $\mathcal{D}^s(S^1)/S^1$, 
i.e. `modulo rotations' and hence one considers all solutions modulo 
traveling.} 

In Section \ref{curvature} we compute the sectional curvature of 
$\mathcal{D}^s(S^1)$ and in Section \ref{Virasoro} we revisit 
the derivation of the $\mu$HS equation on the Virasoro group and 
describe its two compatible Poisson structures. 

Finally, we remark that although our focus in this paper is on 
the periodic case many of the arguments presented here can be 
adapted, with obvious modifications, to study the corresponding 
equation on the real line. 

\bigskip


\section{Eulerian nature of the $\mu$HS equation}\label{GeoSec}\nequation

\subsection{Geodesic flow on the group of circle diffeomorphisms}
Let $\mathcal{D}^s(S^1)$ denote the set of circle diffeomorphisms 
of a sufficiently high Sobolev class $H^s$ and let 
$T\mathcal{D}^s(S^1)$ be its tangent bundle.\footnote{In the bulk of 
the paper we can assume that the Sobolev index $s>3/2$ is simply large 
enough for the computations to make sense 
(for example, taking $s=\infty$, which corresponds to the Fr\'echet 
topology on smooth functions, will do). 
In Section \ref{WP} on well-posedness, where analytical details 
are of essence, we will be more precise.}
Consider the following inner product on the corresponding Lie
algebra of vector fields, that is, on the tangent space
$T_{id} \mathcal{D}^s(S^1)
=
\{u(x)\partial/\partial x:  u \in H^s(S^1)\}$:
\begin{equation}\label{metric}
\langle u, v \rangle = \mu(u) \mu(v) + \int u_x v_x \, dx
=
\int u Av \, dx\,,
\end{equation}
where the inertia operator $A = \mu - \partial_x^2$ acts on a function 
$v$ by $Av = \mu(v) - v_{xx}$.
Using right translations we extend this inner product to 
a right-invariant metric on the group $\mathcal{D}^s(S^1)$ 
so that  
\begin{equation} \label{Metric} 
\langle U, V \rangle_\eta
= 
\int U\eta_x \, dx \int V \eta_x \, dx 
+ 
\int\frac{U_xV_x}{\eta_x}\, dx 
\end{equation} 
for any diffeomorphism $\eta \in \mathcal{D}^s(S^1)$ and 
any $U, V \in T_\eta \mathcal{D}^s(S^1)$. 

\begin{theorem}\label{Eulerth}
The $\mu$HS equation is an Euler equation on $\mathcal{D}^s(S^1)$. 
More precisely, it corresponds to the equation of the geodesic flow 
on $\mathcal{D}^s(S^1)$ with respect to the right-invariant metric 
(\ref{Metric}). 
\end{theorem}

\begin{remark}
{\rm
Observe that $A$ is invertible as an operator from $H^s$ onto $H^{s-2}$ 
with the inverse $A^{-1}w = u$ given explicitly by
\begin{equation} \label{L_inv}
\begin{split} 
u(x) 
=  
\biggl(\frac{x^2}{2} - \frac{x}{2} + \frac{13}{12}\biggr)\mu(w) 
&+ \bigl(x-\frac{1}{2}\bigr) \int\int_0^y w(s) \, ds dy	  \\
&- \int_0^x\int_0^y w(s) \, ds dy 
+ 
\int\int_0^y \int_0^s w(r) \, drdsdy.
\end{split}
\end{equation}
Moreover, note that $A^{-1}$ commutes with $\partial_x$ and 
$\mu(u) = \mu(Au)$. 
}
\end{remark} 

\proofbegin of Theorem \ref{Eulerth}.
Using (\ref{metric}) and the definition of the coadjoint operator 
$$
\langle ad^\ast_v(u), w \rangle = \langle u, [v,w]\rangle,
$$
where $[v, w] = v_xw - vw_x$, we readily compute 
\begin{equation}\label{Bdef}  
ad^\ast_v(u) = A^{-1}(2v_xAu + vAu_x). 
\end{equation}
It follows that the corresponding Euler equation is 
$$
u_t =-ad^\ast_u(u) =-A^{-1}(2\mu(u) u_x -2u_x u_{xx} -u u_{xxx}),
$$
or, applying $A$ to both sides, equivalently 
$$
\mu(u_t) - u_{txx} = - 2 \mu(u) u_x + 2u_x u_{xx} + u u_{xxx}.
$$
Integrating both sides of this equation over the circle and using 
periodicity implies that $\mu(u_t) = 0$, thus yielding 
the $\mu$HS equation. 
\proofend 

An alternative derivation proceeds by calculating the first 
variation of the energy functional of the right-invariant 
metric (\ref{Metric}) 
$$ 
E(\xi) 
= 
\frac{1}{2} \int_0^1 \langle \dot\xi(t),\dot\xi(t)\rangle dt 
$$
defined on the space of smooth paths $t \to \xi(t)$ between 
$\xi(0)=\mathrm{id}$ and $\xi(1)$ in the diffeomorphism group. 
The resulting geodesic equation in $\mathcal{D}^s(S^1)$ 
right-translated to the tangent space at the identity is 
the $\mu$HS equation. 

Combining Theorem \ref{Eulerth} with a local well-posedness result 
of Section \ref{WP} we will show that the metric (\ref{Metric}) 
has a well-defined Riemannian exponential map on 
$\mathcal{D}^s(S^1)$ whenever $s>3/2$ 
(Corollary \ref{exp}).


\subsection{$\mu$HS as an evolution equation}\label{muHSasevolutionsubsec}

Recall that the periodic Hunter-Saxton equation 
$$
u_{txx} + 2u_x u_{xx} + u u_{xxx} = 0
$$
is defined only as an evolution equation on the quotient 
$\mathcal{D}^s(S^1)/S^1$, i.e. 
on the space of cosets 
$\{u(x)+c~|~x\in S^1, \, c \in \R\}$, 
see e.g. \cite{K-M}.
The ambiguity introduced by arbitrariness of $c$ disappears 
in the case of the $\mu$HS equation. 
Indeed, integrating equation (\ref{muHS}) in $x$ gives
\begin{equation}\label{muHSintegrated}
u_{tx} = 2 \mu(u) u  - u u_{xx} - \frac{1}{2}u_x^2  + r_0(t),
\end{equation}
where $r_0(t)$ is uniquely defined at each time $t$ by 
the condition that the right-hand side be a total $x$-derivative. 
Explicitly, one has $r_0(t)= -2\mu(u)^2-\frac 12\mu(u_x^2)$. Note that 
$r_0(t)=r_0$ is constant in time, since 
$-2 r_0(t)=  4\mu(u)^2+\mu(u_x^2)=3\mu(u)^2+\langle u,u\rangle$, while both 
the $\mu$-value and the energy are conserved along $\mu$HS solutions.
Integrating once again over $[0,x]$ gives 
$$
u_{t}(t,x) 
= 
\int_0^x \left(2 \mu(u) u  - u u_{xx} - \frac{1}{2}u_x^2\right)\,dx 
+ 
r_0 x + q_0(t). 
$$ 
The same procedure in the case of the HS equation leads to 
the arbitrary integration constant $q_0(t)$, which is 
the reason for HS being defined only on cosets. 
This time, the condition $\mu(u_t)=0$ 
implies that the mean value of the expression on the right-hand side 
above must be zero. This defines $q_0(t)$ uniquely and hence gives 
a well-defined evolution on the entire space of periodic functions 
$u:S^1 \to \R$.


\section{Variational principles and physical motivation}
\label{varphyssec}\nequation

\subsection{$\mu$HS as a bi-variational equation}
\begin{theorem}
The $\mu$HS equation satisfies two different variational principles.
\end{theorem}

We would like to show that ($\mu$HS) arises as the  equation 
$\delta S=0$ for the action functional 
$$
\mathcal S =\int \big(\int \mathcal L \,dx\big)\,dt
$$
with two different densities $\mathcal L$.

Consider the {\it first Lagrangian density} 
\begin{equation}\label{firstL}
\mathcal L = u_x u_t+uu_x^2+ 2\mu(u)u^2\,.
\end{equation}
This is a generalization of the HS Lagrangian density where the $\mu$-term
is added. The variational principle $\delta \mathcal S=0$
gives the Euler-Lagrange equation
\begin{equation}\label{firstEL}
u_{tx} = - \frac{1}{2}u_x^2 - uu_{xx} + 2\mu(u)u + \mu(u^2).
\end{equation}
The right-hand-side of this equation has to be understood up to addition of 
a (possibly time-dependent) constant to make it a total derivative:
\begin{equation}\label{firstELc}
u_{tx} = - \frac{1}{2}u_x^2 - uu_{xx} + 2\mu(u)u + \mu(u^2)+c\,.
\end{equation}
Here $c$ is defined uniquely  by the condition
$$
\int u_x^2 dx + 4\mu(u)^2 + 2\mu(u^2) +2c = 0\,,
$$
vanishing of the mean of the right-hand-side in (\ref{firstELc}). This
immediately leads to equation (\ref{muHSintegrated}). 
Equivalently, by taking the $x$-derivative,
equation (\ref{firstEL}) can be directly understood as the $\mu$HS 
equation for $ u_{txx}$.

\begin{remark}
{\rm
The additional constant term $c$ in equation (\ref{firstELc}) 
is necessitated by the term $u_x u_t$ in the Lagrangian.
Indeed, consider any variational principle
$$
\delta\left( \int \int u_x u_t dx dt + f(u) \right)= 0.
$$
It gives the equation $u_{tx}=-\frac{1}{2}\delta f/\delta u$, which makes 
sense on the circle only if the mean of the function $\delta f/\delta u$ 
is zero. 
Thus by including the term $u_x u_t$ 
in the Lagrangian, we are {\it assuming} that the Euler-Lagrange equation 
will describe the evolution of the complete derivative $u_x$ 
in the form
$$
u_{tx} = X(u),
$$
for some vector field $X$. By adding an appropriate $c$ for each moment 
of time to the right-hand-side, one can {\it force} the result to be 
a derivative.

For the Lagrangian (\ref{firstL})
this constant will automatically give the appropriate value of $r_0$ 
already discussed in Section \ref{muHSasevolutionsubsec}.
Most importantly, this constant  is non-essential and will disappear once 
we differentiate both sides to obtain the equation for $u_{txx}$, 
as in the original $\mu$HS equation.\footnote{It is worth
mentioning that the same kind of modification by a constant
is required in the variational principle for the original HS equation 
in the periodic case, even without the $\mu$-term.}
}
\end{remark}

\smallskip

The {\it second variational representation} can be obtained from
the Lagrangian density
\begin{equation}\label{secondL}
\bar {\mathcal L} = \frac{1}{2}(u_x^2 + \mu(u) u) + \pi(z_t + uz_x)\,,
\end{equation}
which is a $\mu$-generalization of the density 
presented in \cite{H-Z}. Varying the corresponding action with respect 
to $u, \pi$, respectively $z$, we find
\begin{align}
-u_{xx} + \mu(u) + \pi z_x &= 0,	\label{first}
		\\ 
z_t + uz_x &= 0,			\label{second}
		\\		
-\pi_t - (\pi u)_x &= 0.		\label{third}
\end{align}
Thus, using (\ref{first}), the expression $u_{txx} + 2u_xu_{xx} + uu_{xxx}$ 
equals
$$
\pi_tz_x + \pi z_{tx} + 2u_x(\mu(u) + \pi z_x) + u(\pi z_x)_x.
$$
After employing (\ref{second}) and (\ref{third}) to replace $z_{tx}$ 
respectively $\pi_t$, this simplifies to just $2\mu(u) u_x$, showing that 
the equations of motion is the $\mu$HS equation.

As we will see in Section \ref{bihambivarsec}, these bivariational principles 
correspond to the two Hamiltonian formulations of the $\mu$HS equation.

\medskip


\subsection{Physical motivation}
First recall the derivation of the HS equation from \cite{H-S}.
One considers the director field $\mathbf{n}(x,y,z,t)$ which describes the 
average orientation of the long rod-like molecules in a nematic liquid crystal. 
Assuming invariance under spatial translations and no fluid flow, the action 
is given by (see (A.14) in \cite{H-S})
$$\mathcal{S} 
= \int \int \biggl\{ \frac{1}{2}|\mathbf{n}_t|^2 
- W(\mathbf{n}, \nabla \mathbf{n}) 
+ \frac{1}{2}\lambda |\mathbf{n}|^2\biggr\} d\mathbf{x} dt,
$$
where $\lambda$ is a Lagrangian multiplier enforcing $|\mathbf{n}| = 1$, 
and the Oseen-Frank internal energy $W$ is
$$
W(\mathbf{n}, \nabla \mathbf{n}) 
= \frac{1}{2}\biggl(k_1 (\nabla \cdot \mathbf{n})^2 
+ k_2(\mathbf{n} \cdot \nabla \times \mathbf{n})^2 
+ k_3|\mathbf{n} \times (\nabla \times \mathbf{n})|^2\biggr)\,,
$$
where $k_1, k_2,$ and $k_3$, are nonnegative physical parameters 
describing how the material responds to splay, twist, and bend, respectively.
For the director field
$$
\mathbf{n}(x,y,z,t) = \cos \psi(t,x) \mathbf{e}_x + \sin \psi(t,x) \mathbf{e}_y,
$$
we get
$$
W(\mathbf{n}, \nabla \mathbf{n})
 = \frac{1}{2}\biggl(k_1 \sin^2 \psi + k_3 \cos^2 \psi\biggr)\psi_x^2.
$$
Writing
$$
c^2(\psi) = k_1 \sin^2 \psi + k_3 \cos^2 \psi,
$$
it follows that
$$
W(\mathbf{n}, \nabla \mathbf{n}) = \frac{1}{2}c^2(\psi) \psi_x^2.
$$
Hence, since $|\mathbf{n}_t|^2 = \psi_t^2$, 
this leads to the action
$$
\mathcal{S} = 
\int\int \frac{1}{2}\bigl(\psi_t^2 - c^2(\psi) \psi_x^2\bigr)
\,dx \,dt\,.
$$
(Note that the Lagrange multiplier $\lambda$ does not appear because 
$\mathbf{n}$ is already normalized to have length one.)
The Euler-Lagrange equation, obtained in \cite{H-S}
for the action $\mathcal S$ is
$$
\psi_{tt} - c(\psi)[c(\psi)\psi_x]_x = 0.
$$

At this point we deviate from the HS derivation and 
assume that there is a preferred direction $\psi_0$ of the director field.
This can be related, e.g.,  to an exterior magnetic field, acting on dipoles
and returning them to this direction. Such a returning force corresponds to 
the addition of the potential energy term $-\frac{k}{2}(\psi - \psi_0)^2$ 
to the Lagrangian density 
$\mathcal{L}=\frac{1}{2}(\psi_t^2 - c^2(\psi) \psi_x^2)$.
We also assume that there is retroactive action
of other dipoles on a given director, i.e. that  Hooke's elasticity 
coefficient $k :=\mu(\psi - \psi_0)$ is proportional to the average deviation of 
all dipoles from $\psi_0$: other dipoles act on a given one, and 
the more other dipoles are rotated, the stronger the returning force is on
the given dipole.

Thus the new action is
$$
\mathcal{S}=\frac{1}{2}\int \int \bigl(\psi_t^2 - c^2(\psi) \psi_x^2
-\mu(\psi - \psi_0)(\psi - \psi_0)^2 \bigr)\,dx \,dt\,,
$$
where $\psi_0$ is the direction about which the fluctuations 
occur.\footnote{The added term does not preserve the periodicity of 
the angle $\psi$, but this is inessential since when deriving the 
asymptotic equation we only consider small fluctuations of $\psi$.}
We need to choose $\psi_0$ so that both $c_0 := c(\psi_0)$, 
as well as $c'_0 := c'(\psi_0)$, are nonzero. 

The Euler-Lagrange equation now becomes
\begin{equation}\label{modifiedEL}  
  \psi_{tt} - c(\psi)[c(\psi)\psi_x]_x + \mu(\psi - \psi_0)(\psi - \psi_0)
+ \frac{1}{2}\mu( (\psi - \psi_0)^2 ) = 0.
\end{equation}
Now, repeating the derivation in \cite{H-S}
with this modified equation, we introduce the variables 
$$\theta = x - c_0 t, \qquad \tau = \epsilon t,$$
and substitute the expansions
$$
\psi = \psi_0 + \epsilon \psi_1(\theta, \tau) 
+ \epsilon^2 \psi_2(\theta, \tau) + O(\epsilon^3) \quad \hbox{and} \quad
c(\psi) = c_0 + \epsilon c'_0 \psi_1 + O(\epsilon^2)
$$
into (\ref{modifiedEL}). One finds
$$
\bigl(-c_0 \partial_\theta + \epsilon\partial_\tau\bigr)^2\bigl(\epsilon\psi_{1} 
+ \epsilon^2 \psi_2\bigr)
- \bigl(c_0 + \epsilon c'_0 \psi_1\bigr)
\biggl[\bigl(c_0 + \epsilon c'_0 \psi_1\bigr)\bigl(\epsilon\psi_{1\theta} 
+ \epsilon^2 \psi_{2\theta}\bigr)\biggr]_\theta 
$$
$$
+\bigl(\epsilon \psi_1 + \epsilon^2 \psi_2\bigr) \int \bigl(\epsilon \psi_1 
+ \epsilon^2 \psi_2\bigr) dx + 
\frac{1}{2}\int \bigl(\epsilon \psi_1 + \epsilon^2 \psi_2\bigr)^2 dx 
+ O(\epsilon^3) = 0,
$$
where we used that
$$
\partial_t = -c_0 \partial_\theta + \epsilon\partial_\tau
\qquad\text{ and }\qquad
\partial_x = \partial_\theta\,.
$$
Collecting terms of different order we see that both the $O(1)$ and the 
$O(\epsilon)$ levels vanish, whereas the $O(\epsilon^2)$ terms give
$$
c_0^2 \psi_{2\theta \theta} - 2c_0\psi_{1\theta \tau} 
- c_0 c'_0 \bigl[ \psi_1\psi_{1\theta}\bigr]_\theta 
- c_0^2 \psi_{2\theta\theta} - 
c_0 c'_0 \psi_1 \psi_{1\theta\theta}
+ \psi_1 \int \psi_1 dx + \frac{1}{2}\int  \psi_1^2 dx = 0.
$$
This is simplified to
$$
\psi_{1\theta \tau} 
= - \frac{c'_0}{2}(\psi_{1\theta}^2 + 2\psi_1 \psi_{1\theta \theta}) 
+ \frac{\psi_1}{2c_0} \int \psi_1 dx + \frac{1}{4c_0} \int  \psi_1^2 dx.
$$
Replacing $\tau \to t$, $\theta \to x$, and $\psi_1 \to u/c'_0$, we arrive at
$$
u_{tx} = - \frac{1}{2}u_x^2 - uu_{xx} + \frac{1}{2c_0c'_0} (u \mu( u)
+ \frac{1}{2}\mu( u^2)).
$$
This is the $\mu$HS equation in the form (\ref{firstEL}) 
when $\frac{1}{2c_0c'_0} = 2$.
After differentiation this leads to the 
evolution equation ($\mu$HS) for $u_{txx}$.

\bigskip


\section{Hamiltonian framework and conservation laws} 
\label{Biham} \nequation


\subsection{Bihamiltonian structure and the $\mu$HS hierarchy}
\label{manyham}

\begin{theorem} The $\mu$HS equation is Hamiltonian with respect to two
compatible Hamiltonian structures.
\end{theorem}

\proofbegin ~ 
Let $m := Au= \mu(u) - u_{xx}$ and define
$$H_1 = \frac{1}{2}\int um dx = \frac{1}{2}\langle u, u\rangle, \qquad 
H_2 = \int (\mu(u)u^2 + \frac{1}{2}uu_x^2 )dx.$$
Then
$$
\frac{\delta H_1}{\delta u} = m,\qquad \frac{\delta H_2}{\delta u} 
= \mu(u^2) + 2\mu(u)u - \frac{1}{2}u_x^2 - u u_{xx}.$$
Since $A$ is self-adjoint, we have the relation 
$\frac{\delta H_1}{\delta u} = A\frac{\delta H_1}{\delta m}$. Hence
\begin{equation}\label{variationalderivatives}
  \frac{\delta H_1}{\delta m} = u,\qquad \frac{\delta H_2}{\delta m} 
= A^{-1}(\mu(u^2) + 2\mu(u)u - \frac{1}{2}u_x^2 - u u_{xx}),
\end{equation}
so that the $\mu$HS equation
$$
-u_{txx} = -2\mu(u)u_x + 2u_xu_{xx} + uu_{xxx},
$$
can be written both as
$$
m_t = \mathcal{B}_1\frac{\delta H_1}{\delta m},
$$
and as
$$
m_t = \mathcal{B}_2\frac{\delta H_2}{\delta m},
$$
where the Hamiltonian operators are
$$
\mathcal{B}_1 = -(m\partial_x + \partial_x m), 
\qquad \mathcal{B}_2 = -\partial_xA = \partial_x^3.
$$
These operators are compatible as a straightforward but lengthy calculation shows. 
We will present their geometric meaning in the more general setting of the Virasoro
algebra in Section \ref{Virasoro}.
\proofend

\medskip

Using this bihamiltonian structure it is well-known that one can 
(at least formally) obtain an infinite hierarchy of conservation laws, 
cf. \cite{McK2}. 
If two functionals $F_0$ and $F_1$ satisfy
$$
\mathcal{B}_1\frac{\delta F_0}{\delta m} 
= \mathcal{B}_2\frac{\delta F_1}{\delta m},
$$
we say that $F_0$ raises to $F_1$ or that $F_1$ lowers to $F_0$; in symbols 
$F_0 \nearrow F_1$ or $F_1 \searrow F_0$.
Starting from $H_1$ and $H_2$ the goal is to construct an infinite
 sequence of conservation laws such that
$$
\cdots \nearrow H_{-1} \nearrow H_0 \nearrow H_1 \nearrow H_2 
\nearrow H_3 \nearrow \cdots.
$$

Just like for the Camassa-Holm equation, we find that lowering of 
the $H_n$'s is unobstructed for the $\mu$HS equation, i.e.
there are $H_{-1}, H_{-2}, \dots$, such that $H_{-1} \searrow H_{-2}
\searrow \dots$. Moreover, $H_{-1}, H_{-2}, \dots$ are local functionals
given by
\begin{equation}\label{H-n}
   H_{-n} 
= \frac{1}{3/2-n} \int m \frac{\delta H_{-n}}{\delta m} dx, \qquad n = 1,2,\dots,
\end{equation}
where $\frac{\delta H_{-n-1}}{\delta m}$ is obtained from 
$\frac{\delta H_{-n}}{\delta m}$ via
$$
\frac{\delta H_{-n-1}}{\delta m} 
= \mathcal{B}_1^{-1}\mathcal{B}_2\frac{\delta H_{-n}}{\delta m} 
= -\frac{1}{2\sqrt{m}}\int^x \frac{1}{\sqrt{m}}
\left(\frac{\delta H_{-n}}{\delta m}\right)_{xxx} dx.
$$
Starting with $\delta H_{-1} = \frac{1}{2 \sqrt{m}}$ this provides 
a constructive approach for finding the $H_{-n}$'s, which is easily implemented in Mathematica. For example, we get
$$
\frac{\delta H_{-2}}{\delta m} 
= \frac{m_{xx}}{8m^{5/2}} - \frac{5m_x^2}{32m^{7/2}},
$$
$$
\frac{\delta H_{-3}}{\delta m} = 
\frac{1155m_x^4}{1024m^{13/2}} - \frac{231m_x^2m_{xx}}{128m^{11/2}} +
\frac{21m_{xx}^2}{64m^{9/2}} + \frac{7m_xm_{xxx}}{16m^{9/2}} -  
\frac{m_{xxxx}}{16m^{7/2}},$$
and the corresponding functionals $H_{-2}$ and $H_{-3}$ 
are obtained from (\ref{H-n}).

The situation for raising is more intricate due to the occurence of nonlocal 
functionals. Nevertheless, we can deduce the existence of expressions $Q_n[m]$ 
(possibly nonlocal in $m$) such that $Q_0 = \frac{\delta H_{0}}{\delta m}$ and
$$
Q_{n+1} = \mathcal{B}_2^{-1}\mathcal{B}_1Q_n, \qquad n \geq 0.
$$
If there exist functionals $H_n[m]$ with $\frac{\delta H_n}{\delta m} = Q_n$, 
they constitute an infinite ladder of functionals with 
$H_1 \nearrow H_2 \nearrow H_2 \nearrow \cdots.$ 
The existence of the $H_n$ is not evident in this approach
due to the nonlocal nature of the $Q_n$'s.
However, it follows from the expansion of a Casimir for a family 
of compatible Poisson brackets in the approach of Section \ref{Virasoro}, cf. 
\cite{K-M}.

In summary, the first few conservation laws in the hierarchy are
\begin{eqnarray*}
H_{-2} = -\frac{1}{16}\int \frac{m_x^2}{m^{5/2}} dx,
&& \frac{\delta H_{-2}}{\delta m} = \frac{m_{xx}}{8m^{5/2}} - \frac{5m_x^2}{32m^{7/2}},
     	\\
H_{-1} = \int \sqrt{m} dx,
  && \frac{\delta H_{-1}}{\delta m} = \frac{1}{2 \sqrt{m}},
     	\\
H_0 = \int {m dx},
&& \frac{\delta H_{0}}{\delta m} = 1,
     	\\
H_1 = \frac{1}{2}\int {mu dx},
&& \frac{\delta H_{1}}{\delta m} = u,
     	\\
H_2 = \int \left(\mu(u)u^2 + \frac{1}{2}uu_x^2 \right)dx,
&& \frac{\delta H_{2}}{\delta m} = (\mu - \partial_x^2)^{-1}\left(\mu(u^2) 
+ 2\mu(u)u - \frac{1}{2}u_x^2 - u u_{xx}\right).
\end{eqnarray*}
with corresponding commuting flows
\begin{align*}
m_t =& \mathcal{B}_1\frac{\delta H_{-2}}{\delta m} = \mathcal{B}_2\frac{\delta H_{-1}}{\delta m}
= \biggl(\frac{1}{2\sqrt{m}}\biggr)_{xxx},
		\\
m_t =& \mathcal{B}_1\frac{\delta H_{-1}}{\delta m} = \mathcal{B}_2\frac{\delta H_{0}}{\delta m} = 0,
		\\
m_t =& \mathcal{B}_1\frac{\delta H_{0}}{\delta m} = \mathcal{B}_2\frac{\delta H_{1}}{\delta m} = -m_x,
		\\
m_t =& \mathcal{B}_1\frac{\delta H_{1}}{\delta m} = \mathcal{B}_2\frac{\delta H_{2}}{\delta m} = -2mu_x - m_x u.
\end{align*}


\subsection{Relation with bi-variational principles}\label{bihambivarsec}
Recall that the $\mu$HS equation can be obtained from two distinct 
variational principles, corresponding to the Lagrangian densities 
$\mathcal{L}$ and $\bar{\mathcal{L}}$ introduced in Section \ref{varphyssec}.

\begin{theorem}
The two variational formulations for the $\mu$HS equation formally
correspond to the two Hamiltonian formulations of this equation
with Hamiltonian functionals $H_1$ and $H_2$.
\end{theorem}

\proofbegin ~ 
The action is related to the Lagrangian by 
$S = \int L dt$ and $L = \int \mathcal{L} dx$.
The first variational principle has the Lagrangian density of equation (\ref{firstL}),
$$
\mathcal{L} = u_x u_t+uu_x^2+ 2\mu(u)u^2\,,
$$
where we skip the linear term.
The momenta conjugate to the velocities $u_t$ are
$$
\frac{\partial \mathcal{L}}{\partial u_t} = u_{x}.
$$
Hence the Hamiltonian density is
$$
\mathcal{H} = u_{x}u_t - \mathcal{L} = -uu_x^2- 2\mu(u)u^2,
$$
and so
$$
H = -\int (uu_x^2 + 2\mu(u) u^2) dx,
$$
which indeed agress with the conservation law $H_2$ up to a factor of $-2$. 

\medskip

In the second principle the Lagrangian density is
$$
\bar{\mathcal{L}} = \frac{1}{2}u_x^2 + \frac{1}{2} \mu(u) u + \pi(z_t + uz_x)\,.
$$
The momenta conjugate to the velocities $u_t, \pi_t$, and $z_t$,
respectively, are
$$
\frac{\partial \bar{\mathcal{L}}}{\partial u_t} = 0, \quad
\frac{\partial \bar{\mathcal{L}}}{\partial \pi_t} = 0, \quad
\frac{\partial \bar{\mathcal{L}}}{\partial z_t} = \pi.
$$
Consequently, the Hamiltonian density is
$$
\mathcal{H} = z_t \pi - \bar{\mathcal{L}} 
= z_t \pi - \left(\frac{1}{2}u_x^2 + \frac{1}{2}\mu(u) u + \pi(z_t + uz_x)\right).
$$
Using the equation of motion (\ref{first}) to substitute for $\pi z_x$, 
this becomes
$$
\mathcal{H} = -\frac{1}{2}u_x^2 - uu_{xx} + \frac{1}{2}\mu(u) u.
$$
Therefore, the Hamiltonian is
$$
H = \int \mathcal{H} dx = \frac{1}{2} \int (u_x^2 + \mu(u) u) dx,
$$
which is exactly $H_1$.
\proofend

\bigskip

\begin{remark}(Lax Pair)
{\rm
The $\mu$HS equation
$$
u_{txx} - 2\mu(u) u_x + 2u_x u_{xx} + u u_{xxx} = 0.
$$
can be viewed as the condition of compatibility between
$$\psi_{xx} = \lambda m \psi,$$
and
$$\psi_t = \biggl(\frac{1}{2\lambda} - u\biggr)\psi_x + \frac{1}{2}u_x \psi,$$
where $\lambda$ is a spectral parameter and, as above, $m = \mu(u) - u_{xx}$. 
In other words, assuming that these equations are satisfied for 
all nonzero $\lambda \in \C$, 
the relation $\psi_{txx} = \psi_{xxt}$ holds
if and only if $u(t,x)$ satisfies the $\mu$HS equation.
}
\end{remark}


\section{The periodic Cauchy problem} 
\label{WP} \nequation

In this section we discuss well-posedness of the $\mu$HS equation. 
We do not aim for the strongest possible results here. 
Instead, after establishing a basic local well-posedness theorem 
we proceed to show that while some solutions persist indefinitely 
other solutions develop singularities in finite time. 

Recall that the mean of any solution $u=u(t,x)$ 
is conserved by the flow and hence the Cauchy problem for 
the $\mu$HS equation can be recast in the form 
\begin{equation} \label{LmuHS} 
u_t + uu_x 
+ 
A^{-1}\partial_x \Big( 2\mu(u)u +\frac{1}{2}u_x^2 \Big) 
= 0  
\end{equation} 
and 
\begin{equation} \label{ic} 
u(0,x) = u_o(x), 
\qquad 
x \in S^1, \quad  t \in \mathbb{R}\,,
\end{equation} 
where $A= \mu - \partial_x^2$ is an isomorphism 
between $H^s(S^1)$ and $H^{s-2}(S^1)$ with inverse $u=A^{-1}w$
given explicitly by equation (\ref{L_inv}).
Conservation of the mean will play an important role later on. 


\subsection{Local well-posedness and persistence result} 
In order to establish local well-posedness of the Cauchy problem 
(\ref{LmuHS})--(\ref{ic}) we adopt here the approach developed 
in \cite{Mis02} for the CH equation. 
Our main result is the following 

\begin{theorem} \label{wp_PDE} 
If $s>3/2$ then for any $u_o \in H^s(S^1)$ there exist $T>0$ 
and a unique solution 
$u \in C((-T,T),H^s)\cap C((-T,T),H^{s-1})$ of the initial value problem 
(\ref{LmuHS})--(\ref{ic}) which depends continuously 
on the initial data $u_o$. 
Furthermore, the solution persists as long as 
$\|u(t,\cdot)\|_{C^1}$ stays bounded. 
\end{theorem}

First, however, we will study the associated Cauchy problem 
on the group $\mathcal{D}^s(S^1)$ of Sobolev $H^s$ diffeomorphisms 
with $s>3/2$. 
Given a solution $u(t,x)$ of (\ref{LmuHS}) with initial data $u_o$ 
we let $t \to \eta(t,x)$ be the flow of $u$, that is 
\begin{equation} \label{flow}
\dot{\eta}(t,x) 
= 
u(t,\eta(t,x)). 
\end{equation} 
Differentiating both sides of (\ref{flow}) with respect to $t$ 
we rewrite the original Cauchy problem (\ref{LmuHS})--(\ref{ic}) as 
\begin{equation} \label{ODE} 
\ddot{\eta}
= 
- \left\{ A^{-1} \partial_x 
\Big( 2\mu(\dot\eta\circ\eta^{-1})\,\dot{\eta}\circ\eta^{-1} 
+
\frac{1}{2}\big( 
(\dot{\eta}\circ\eta^{-1})_x 
\big)^2 \Big) \right\} \circ\eta 
=: f(\eta, \dot{\eta}) 
\end{equation} 
\begin{equation}  \label{ODEic} 
\eta(0,x) = x, 
\quad 
\dot\eta(0,x) = u_o(x) 
\end{equation} 
and prove 

\begin{theorem} \label{wp_ODE} 
Let $s>3/2$. For any $u_o \in H^s(S^1)$ there exist $T>0$ and 
a unique solution 
$(\eta, \dot{\eta}) \in C^1((-T,T),\mathcal{D}^s(S^1)\times H^s)$ 
of the initial value problem (\ref{ODE})--(\ref{ODEic}) which depends 
differentiably on the initial data $u_o$. 
\end{theorem} 

An immediate consequence of Theorem \ref{wp_ODE} combined with 
the inverse function theorem for Banach spaces 
is the existence of the Riemannian exponential map in 
$\mathcal{D}^s(S^1)$ whose geodesics were described 
in Theorem \ref{Eulerth}. 
The corollary below permits a rigorous development of 
the associated infinite-dimensional Riemannian geometry and justifies 
the geometric constructions of Sections 2 and 7. 

\begin{corollary} \label{exp} 
For any $s>3/2$ the metric (\ref{Metric}) has a well-defined 
exponential map 
$
\exp_{\mathrm{id}}: 
U \subset T_{\mathrm{id}}\mathcal{D}^s(S^1) \to \mathcal{D}^s(S^1)
$
which is a $C^1$ diffeomorphism from an open set $U$ onto 
a neighbourhood of the identity in $\mathcal{D}^s(S^1)$. 
\end{corollary}

Our first main task will be to show that the map 
$(\xi, w) \to \big( w, f(\xi,w) \big)$ 
defines a $C^1$ vector field on the Hilbert manifold 
$\mathcal{D}^s(S^1)\times H^s(S^1)$. 
The proof of this fact is similar to the corresponding proof for 
CH and we will refer the reader to \cite{Mis02} for some of the details. 

\medskip
\proofbegin of Theorem \ref{wp_ODE}. 
It is sufficient to show that $f$ is continuously differentiable 
(in the Fr\'echet sense) with respect to the variables $\xi$ and $w$ 
in some neighbourhood of the point $(\mathrm{id}, 0)$ 
in $\mathcal{D}^s(S^1)\times H^s(S^1)$. 
The following technical result will be useful. 

\begin{lemma}\label{H_lemma} {\rm (\cite{Mis02}, Appendix 1)}
For $s>3/2$ the composition map $\xi \to w\circ\xi$ with 
an $H^s$ function $w$ and the inversion map $\xi \to \xi^{-1}$ 
are continuous from $\mathcal{D}^s(S^1)$ to $H^{s}(S^1)$ and 
from $\mathcal{D}^s(S^1)$ to itself respectively.
Moreover 
\begin{equation} \label{lem} 
\|w\circ\xi\|_{H^s} 
\leq 
C \big(1 + \|\xi\|_{H^s}^s\big) \|w\|_{H^s} 
\end{equation} 
with $C$ depending only on $\inf{|\partial_x\xi|}$ 
and $\sup{|\partial_x\xi|}$. 
\end{lemma} 

In particular, using (\ref{lem}) we find that the right hand side 
of (\ref{ODE}) maps $(\xi,w) \to f(\xi,w)$ into $H^s(S^1)$. 

Next, rewrite $f$ in the form 
\begin{equation} \label{f} 
f(\xi,w) = 
- A_{\xi}^{-1} \, \partial_{x,\xi} \, h(\xi,w) \,,
\end{equation}
where 
$
\partial_{x,\xi} = dR_\xi \circ \partial_x \circ dR_{\xi^{-1}} 
$
and 
$
A_{\xi}^{-1} = dR_\xi\circ A^{-1}\circ dR_{\xi^{-1}} 
$ 
are the conjugations of $\partial_x$ and $A^{-1}$ 
by a diffeomorphism $\xi$ 
(here $R_\xi$ denotes the right translation by $\xi$ in 
$\mathcal{D}^s(S^1)$) and where 
\begin{equation} \label{h}
h(\xi,w) 
= 
2 w \int w\circ\xi^{-1}dx 
+ 
\frac{1}{2} \left( \frac{\partial_x w}{\partial_x\xi} \right)^{2}. 
\end{equation}
Let $s \to \xi_s$ be a smooth curve in $\mathcal{D}^s(S^1)$ 
passing through the point $\xi_0 = \xi$ in the direction 
$\partial_{s}\xi_{s}\vert_{s=0}=v$. 
A formal computation gives 
$$
\partial_\xi(\partial_{x,\xi})(v) 
=
\partial_s (\partial_{x,\xi_s})\vert_{s=0} 
= 
\big[ v\circ\xi^{-1} \partial_x , \partial_x \big]_\xi 
$$
and 
$$
\partial_\xi(A^{-1}_{\xi})(v) 
= 
\partial_s (A^{-1}_{\xi_s})\vert_{s=0} 
= 
- A^{-1}_{\xi} \, 
\partial_s (A_{\xi_s})\vert_{s=0} \, 
A^{-1}_{\xi} 
= 
-A^{-1}_{\xi} 
\big[ v\circ\xi^{-1}\partial_x, A \big]_\xi 
A^{-1}_{\xi} \,,
$$
where $[\cdot , \cdot ]$ denotes the commutator of operators. 
Similarly from (\ref{h}) we find 
\begin{equation} \label{partial_xi_h}
\begin{split} 
\partial_\xi h_{(\xi,w)}(v) 
= 
2\Big\{ 
w\circ\xi^{-1} \int w\circ\xi^{-1}\partial_x(v\circ\xi^{-1})dx 
- 
\frac{1}{2} \big(\partial_x(w\circ\xi^{-1})\big)^2 
\partial_x(v\circ\xi^{-1}) 
\Big\}\circ\xi \,,
\end{split}
\end{equation}
as well as 
\begin{equation} \label{partial_w_h}
\begin{split} 
\partial_w h_{(\xi,w)}(v) 
=
2 \Big\{ 
w\circ\xi^{-1} \int v\circ\xi^{-1} dx 
&+ 
v\circ\xi^{-1}\int w\circ\xi^{-1} dx    \\ 
&+ 
\frac{1}{2}\partial_x(w\circ\xi^{-1}) \partial_x(v\circ\xi^{-1}) 
\Big\} \circ\xi. 
\end{split}
\end{equation}
With the help of the above formulas we can compute the directional 
derivatives of $f$ in (\ref{f}), namely 
\begin{equation} \label{partial_xi_f} 
\begin{split} 
\partial_\xi f_{(\xi,w)} (v) 
&= 
- \partial_\xi (A^{-1}_{\xi})(v) 
\partial_{x,\xi} h(\xi,w) 
- 
A^{-1}_{\xi} \partial_\xi(\partial_{x,\xi})(v) h(\xi,w) 
- 
A^{-1}_{\xi}\partial_{x,\xi} \partial_\xi h_{(\xi,w)}(v)  \\ 
&=
\Big\{ 
A^{-1}\partial_x\big( 
v\circ\xi^{-1}\partial_x(h\circ\xi^{-1}) \big) 
- 
v\circ\xi^{-1}\partial_x A^{-1}\partial_x(h\circ\xi^{-1}) 
\Big\}\circ\xi                 \\ 
- 
2\Big\{A^{-1} &\partial_x \Big( 
w\circ\xi^{-1} \int w\circ\xi^{-1}\partial_x(v\circ\xi^{-1})dx 
- 
\frac{1}{2}(\partial_x(w\circ\xi^{-1}))^2 
\partial_x(v\circ\xi^{-1}) 
\Big) \Big\}\circ\xi
\end{split} 
\end{equation} 
and 
\begin{equation} \label{partial_w_f} 
\begin{split} 
\partial_w f_{(\xi,w)}(v) 
&= 
-\big( A^{-1} \partial_x \big)_\xi 
\partial_w h_{(\xi,w)}(v)      \\ 
&\begin{split} 
= 
-2\Big\{A^{-1} \partial_x \Big( 
w\circ\xi^{-1} \int v\circ\xi^{-1} dx 
&+ 
v\circ\xi^{-1}\int w\circ\xi^{-1} dx    \\ 
&+ 
\frac{1}{2}\partial_x(w\circ\xi^{-1}) \partial_x(v\circ\xi^{-1}) 
\Big)\Big\} \circ\xi. 
\end{split}
\end{split} 
\end{equation} 

Our strategy will be to show that both directional derivatives 
$\partial_\xi f_{(\xi,w)}$ and $\partial_w f_{(\xi,w)}$ 
are bounded operators in $L(H^s,H^s)$ that depend continuously 
on $(\xi,w)$ in some neighbourhood of $(\mathrm{id},0)$. 
Using the identities 
\begin{equation} \label{id_1} 
A^{-1}\partial_x w(x) 
= 
\Big( x - \frac{1}{2} \Big) \int w(x) \, dx 
- 
\int_0^x w(y) \, dy + \int\int_0^x w(y) \, dy dx 
\end{equation} 
and 
\begin{equation} \label{id_2} 
A^{-1} \partial_x^2 w(x) 
= 
- w(x) + \int w(x) \, dx \,,
\end{equation} 
which are easily obtained from (\ref{L_inv}) 
and the fact that $A^{-1}$ and $\partial_x$ 
commute, we find explicit expressions for the directional 
derivatives in (\ref{partial_xi_f}) and (\ref{partial_w_f}). 
After a lengthy, but straightforward, calculation 
we obtain 
\begin{equation} \label{partial_xi} 
\begin{split} 
\partial_\xi f_{(\xi,w)}(v)   
=&
\Big\{ 
4\Big( x - \frac{1}{2} \Big) 
\int w\circ\xi^{-1} dx 
\int v\circ\xi^{-1} \partial_x(w\circ\xi^{-1}) \, dx   \\ 
+&
\frac{1}{2}\Big( x -\frac{1}{2} \Big) \int 
\big(\partial_x(w\circ\xi^{-1})\big)^2 \partial_x(v\circ\xi^{-1}) 
\, dx     \; \\
 +  
2\int w\circ\xi^{-1} \partial_x(v\circ\xi^{-1}) \,dx 
&\int_0^x w\circ\xi^{-1} dy     
-  
2\int w\circ\xi^{-1} dx 
\int_0^x v\circ\xi^{-1} \partial_x(w\circ\xi^{-1}) \,dy\\    
- &
\frac{1}{2}\int_0^x 
\big(\partial_x(w\circ\xi^{-1})\big)^2 \partial_x(v\circ\xi^{-1}) 
\, dy         \\ 
+&
2\int w\circ\xi^{-1} dx 
\int\int_0^x v\circ\xi^{-1} \partial_x(w\circ\xi^{-1})\, dy dx  \\
-&
2\int w\circ\xi^{-1} \partial_x(v\circ\xi^{-1})\, dx 
\int\int_0^x w\circ\xi^{-1} dy dx   \\ 
+ 
\frac{1}{2}\int v\circ\xi^{-1} 
\big(\partial_x&(w\circ\xi^{-1})\big)^2 dx 
+
\frac{1}{2}\int\int_0^x 
\big(\partial_x(w\circ\xi^{-1})\big)^2\partial_x(v\circ\xi^{-1})\, dy dx  \\ 
+& 
2 v\circ\xi^{-1} w\circ\xi^{-1} \int w\circ\xi^{-1} dx 
-
2 v\circ\xi^{-1} \Big( \int w\circ\xi^{-1} dx \Big)^2  \\
- &
\frac{1}{2} v\circ\xi^{-1} 
\int \big(\partial_x(w\circ\xi^{-1})\big)^2 dx 
\Big\}\circ\xi
\end{split}
\end{equation} 
and similarly 
\begin{equation} \label{partial_w} 
\begin{split} 
\partial_w f_{(\xi,w)}(v) 
= 
&- \Big\{ 
4\Big(x-\frac{1}{2}\Big) \int w\circ\xi^{-1} \, dx 
\int v\circ\xi^{-1}\, dx   \\ 
&+ 
\Big(x- \frac{1}{2} \Big)
\int \partial_x (w\circ\xi^{-1}) \partial_x(v\circ\xi^{-1}) dx 
- 
2 \int w\circ\xi^{-1} \, dx \int_0^x v\circ\xi^{-1}\, dy      \\ 
&-
2 \int v\circ\xi^{-1}\, dx \int_0^x w\circ\xi^{-1}\, dy 
-
\int_0^x \partial_x(w\circ\xi^{-1}) \partial_x(v\circ\xi^{-1}) dy   \\ 
+ 
2 \int & w\circ\xi^{-1}\, dx \int \int_0^x v\circ\xi^{-1}\, dy dx 
+ 
2 \int v\circ\xi^{-1}\, dx \int\int_0^x w\circ\xi^{-1}\, dy dx \\
&+ 
\int\int_0^x 
\partial_x(w\circ\xi^{-1}) \partial_x(v\circ\xi^{-1})
\, dy dx     
\Big\}\circ \xi. 
\end{split}
\end{equation} 

By Lemma \ref{H_lemma}, in order to show that 
$v \to \partial_\xi f_{(\xi,w)}(v)$ is a bounded operator on $H^s$ 
it will be sufficient to estimate the sum 
$$
\left\| \partial_\xi f_{(\xi,w)}(v) \circ \xi^{-1}\right\|_{L^2} 
+
\left\| 
\partial_x \left( \partial_\xi f_{(\xi,w)}(v)\right)\circ\xi^{-1}
\right\|_{H^{s-1}}. 
$$

Using Cauchy-Schwarz and the Sobolev lemma we can bound the first 
of the two terms above by 
\begin{equation*} 
\begin{split} 
&\lesssim  
\big(\|w\circ\xi^{-1}\|_\infty \|w\circ\xi^{-1}\|_{H^1} 
+
\|w\circ\xi^{-1}\|_{H^1}^2 \big) 
\big(\|v\circ\xi^{-1}\|_{H^1} + \|v\circ\xi^{-1}\|_{C^1}\big)   \\ 
&\lesssim 
\|w\circ\xi^{-1}\|_{H^1}^2 \|v\circ\xi^{-1}\|_{H^s}
\end{split}
\end{equation*} 
and, after differentiating (\ref{partial_xi}) in $x$ and using 
the fact that $H^{s-1}$ is an algebra, estimate the second term by 
\begin{equation*} 
\begin{split} 
&\lesssim 
\big\|\partial_x(w\circ\xi^{-1})\big\|_{L^2}^2
\left\|\partial_x (v\circ\xi^{-1})\right\|_{\infty} 
+ 
\left\|\partial_x(w\circ\xi^{-1})\right\|_{H^{s-1}}^2 
\left\|\partial_x(v\circ\xi^{-1})\right\|_{H^{s-1}}        \\ 
&~~~~~+
\|w\circ\xi^{-1}\|_{\infty}\|w\circ\xi^{-1}\|_{H^{s}}
\|v\circ\xi^{-1}\|_{H^{s}} 
+ 
\|w\circ\xi^{-1}\|_{L^2}^2
\|\partial_x(v\circ\xi^{-1})\|_{H^{s-1}}   \\ 
&\lesssim 
\|w\circ\xi^{-1}\|_{H^s}^2 \|v\circ\xi^{-1}\|_{H^s}. 
\end{split}
\end{equation*} 
Combining these and using the estimate (\ref{lem}) of 
Lemma \ref{H_lemma} we obtain 
$$
\left\|\partial_\xi f_{(\xi,w)}(v)\right\|_{H^s} 
\lesssim 
C \|w\|_{H^s}^2 \|v\|_{H^s} \,,
$$
where the constant $C$ depends only on the $H^s$ norms of 
$\partial_x\xi$ and $\partial_x\xi^{-1}$. 
An analogous argument yields a similar bound on $\partial_w f_{(\xi,w)}$ 
in (\ref{partial_w}) 
$$
\left\|\partial_w f_{(\xi,w)}(v)\right\|_{H^s} 
\lesssim 
C \|w\|_{H^s} \|v\|_{H^s}. 
$$ 
We have thus shown that $f$ is Gateaux differentiable near 
$(\mathrm{id},0)$. 

We next turn to continuity of the directional derivatives 
with respect to the variables $\xi$ and $w$. 
Continuity in $w$ follows from the fact that the dependence 
of both partials on this variable is polynomial. 
Therefore, it only remains to show that the norm of the difference 
$\|\partial_\xi f_{(\xi,w)}(v) 
- 
\partial_\xi f_{(\mathrm{id},w)}(v)\|_{H^s}$ 
can be made arbitrarily small whenever $\xi$ is close to 
$\mathrm{id}$ in $\mathcal{D}^s(S^1)$, uniformly in  
$v, w \in H^s(S^1)$.  
Thus, our task is to estimate the sum 
$$
\left\| \partial_\xi f_{(\xi,w)}(v) 
- 
\partial_\xi f_{(\mathrm{id},w)}(v) \right\|_{L^2}
+ 
\left\| \partial_x\big( \partial_\xi f_{(\xi,w)}(v) 
- 
\partial_\xi f_{(\mathrm{id},w)}(v) \big)\right\|_{H^{s-1}}. 
$$
Using the explicit formulas in (\ref{partial_xi}) together with 
the estimate of Lemma \ref{H_lemma}, the algebra property of $H^{s-1}$ 
and the fact that for any $r>0$ we have a Sobolev embedding into 
H\"older spaces 
$H^r(S^1) \hookrightarrow C^{r-1/2}(S^1)$ with a pointwise estimate 
$$
|g(x) - g(y)| \lesssim \|g\|_{H^r} |x-y|^{r-1/2} 
\qquad\qquad 
( x,y \in S^1 )
$$
both terms in the expression above can be bounded by 
$$
\lesssim 
C\|w\|_{H^s}^2 \|v\|_{H^s} 
\left( 
\big\| \xi^{-1} -\mathrm{id} \big\|_{H^s} 
+ 
\big\| \xi^{-1} -\mathrm{id} \big\|_\infty^{s-3/2}  \right) \,,
$$
where $C$ depends only on the $H^s$ norms of $\xi$ and its inverse. 
The arguments needed to prove continuity of the map 
$\xi \to \partial_w f_{(\xi,w)}$ are very similar. 

From these estimates we conclude that $f$ defined by the right side 
of (\ref{ODE}) is continuously differentiable in a neighbourhood 
of the point $(\mathrm{id},0)$. 
Applying the fundamental ODE theorem for Banach spaces 
we obtain Theorem \ref{wp_ODE}.
\proofend

\medskip

We are now ready to prove local well-posedness 
and persistence results for the $\mu$HS equation. 

\medskip

\proofbegin of Theorem \ref{wp_PDE}.
Observe that directly from the flow equation (\ref{flow}) we have 
$u = \dot\eta \circ\eta^{-1}$. 
Thererfore, since $\mathcal{D}^s(S^1)$ is a topological group 
by Lemma \ref{H_lemma}, the first statement is a direct corollary of 
Theorem \ref{wp_ODE}. 

To prove the persistence result we employ Friedrichs' mollifiers 
$J_\varepsilon \in OPS^{-\infty}$, where $0<\varepsilon <1$, and 
using (\ref{LmuHS}) we first derive an estimate for 
\begin{equation*}  
\begin{split} 
\frac{d}{dt} \|J_\varepsilon u \|_{H^s}^2 
= 
&- 
2\langle 
\Lambda^s J_\varepsilon (u\partial_x u), \Lambda^s J_\varepsilon u 
\rangle_{L^2}  
-
4 \left\langle 
\Lambda^s A^{-1} \partial_x J_\varepsilon 
\big( \mu(u)u \big),\Lambda^s J_\varepsilon u 
\right\rangle_{L^2}    \\ 
&- 
\left\langle \Lambda^s  A^{-1} \partial_x 
J_\varepsilon (\partial_x u)^2, \Lambda^s J_\varepsilon u 
\right\rangle_{L^2},  
\end{split}
\end{equation*}
where $\Lambda^s = \left( 1 - \partial_x^2 \right)^{s/2}$. 
The first term on the right hand side can be bounded as in \cite{tay}, 
see p. 115-116 
$$
\langle 
\Lambda^s J_\varepsilon (u\partial_x u), \Lambda^s J_\varepsilon u 
\rangle_{L^2} 
\lesssim 
\|u\|_{C^1} \|u\|_{H^s}^2. 
$$
For the middle term we have 
$$
\left\langle 
\Lambda^s A^{-1} \partial_x J_\varepsilon 
\big( \mu(u)u \big),\Lambda^s J_\varepsilon u 
\right\rangle_{L^2} 
\lesssim 
\left|\mu(u)\right| 
\left\|\Lambda^s A^{-1}\partial_x J_\varepsilon u\right\|_{L^2} 
\left\| \Lambda^s J_\varepsilon u \right\|_{L^2} 
\lesssim 
\|u\|_{L^1} \|u\|_{H^s}^2. 
$$ 
To estimate the last term we again use Cauchy-Schwarz 
$$
\left\langle \Lambda^s  A^{-1} \partial_x 
J_\varepsilon (\partial_x u)^2, \Lambda^s J_\varepsilon u 
\right\rangle_{L^2} 
\lesssim 
\left\| A^{-1}\partial_x (\partial_x u)^2\right\|_{H^s} 
\|u\|_{H^s}
$$
and with the help of the identities (\ref{id_1}) and (\ref{id_2}) 
we obtain 
\begin{equation*} 
\begin{split} 
\left\| A^{-1}\partial_x (\partial_x u)^2\right\|_{H^s} 
&\simeq 
\left\|A^{-1}\partial_x (\partial_x u)^2\right\|_{L^2} 
+ 
\left\|
\partial_x  A^{-1}\partial_x (\partial_x u)^2 
\right\|_{H^{s-1}}    \\  
&\lesssim 
\|\partial_x u\|_{L^2}^2 
+ 
\big\|-(\partial_x u)^2 +\int (\partial_x u)^2 dx\big\|_{H^{s-1}} \\ 
&\lesssim 
\|\partial_x u\|_{L^2}^2 
+ 
\|\partial_x u \|_\infty \|u\|_{H^s} \,,
\end{split} 
\end{equation*} 
where in the last step we used a standard Moser inequality for products 
of $H^{s-1}$ functions. 

Collecting these estimates and passing to the limit with 
$\varepsilon \to 0$ we arrive at the differential inequality 
$$
\frac{d}{dt} \left\| u \right\|_{H^s}^2 
\lesssim 
\|u\|_{C^1} \|u\|_{H^s}^2\,,
$$ 
which combined with Gronwall's inequality implies that the $H^s$ norm 
of the solution $u$ cannot blow up unless $\|u(t,\cdot)\|_{C^1}$ does. 
\proofend


\subsection{Global solutions with non-negative angular momentum density} 
Our next result shows the effect of the term $2\mu(u)u_x$ on 
the lifespan of solutions to the Cauchy problem. 
\begin{theorem} \label{wp_global} 
Let $s>5/2$. Assume that the function $u_o \in H^s(S^1)$ has 
non-zero mean and satisfies the condition 
\begin{equation} \label{pos} 
A u_o(x) := (\mu - \partial_x^2) u_o(x) \geq 0, 
\qquad 
x \in S^1.
\end{equation} 
Then the Cauchy problem (\ref{LmuHS})--(\ref{ic}) 
has a unique global solution $u$ in the space 
$
C(\mathbb{R}, H^s)\cap C^1(\mathbb{R}, H^{s-1}). 
$
\end{theorem} 

\proofbegin ~ 
In Theorem \ref{wp_PDE} for any such $u_o$ we constructed 
a local solution $u$ with desired regularity and defined up 
to some time $T>0$. 
Furthermore, we showed that in order to extend it beyond $T$ 
it suffices to find a constant $K>0$ such that 
$\|u(t,\cdot)\|_{C^1} \leq K$ for all $|t|<T$. 

Set $m(t,x) = A u(t,x)$ and observe that from 
the identity (\ref{id_1}) we have 
$$ 
\partial_x u(t,x) 
= 
\Big(x - \frac{1}{2}\Big) \int m(t,x)\, dx 
+ 
\int\int_0^x m(t,y) \, dy dx - \int_0^x m(t,y) \, dy 
$$ 
for any $x \in [0,1]$, which consequently yields a bound 
\begin{equation}  \label{bound} 
\| \partial_x u(t) \|_{\infty} 
\lesssim 
\|m(t)\|_{L^1}. 
\end{equation} 
On the other hand, it is not difficult to check that at each point 
of the circle the solution $u$ satisfies a local conservation law 
$$ 
A u(t,\eta(t,x)) (\partial_x \eta(t,x))^2 
= 
A u_o(x)\,,
$$ 
where $\eta$ is the geodesic in $\mathcal{D}^s(S^1)$ corresponding 
to $u$. In particular, if (\ref{pos}) holds then 
$m(t,x) \geq 0$ as long as $u(t,x)$ exists, which in turn yields 
$$
\|m(t)\|_{L^1} 
= 
\int u(t,x) dx = \int u_o(x) \, dx 
$$
and the theorem follows. 
\proofend 

\medskip
An analogous result holds if the assumption (\ref{pos}) 
is replaced with $A u_o(x) \leq 0$. 
Notice that for the periodic Camassa-Holm equation the condition that 
$u-u_{xx}$ does not properly change sign is necessary and sufficient 
for global existence, see \cite{McK}.


\subsection{Breakdown of classical solutions} 
It is not difficult to find solutions which develop 
singularities in finite time. In fact, the conservation of the mean 
\begin{equation} \label{mean} 
\int u(t,x) \, dx = \int u_o(x) \, dx 
\end{equation} 
implies that sufficiently smooth solutions of $\mu$HS corresponding 
to zero-mean initial data must satisfy the periodic HS equation 
\begin{equation*} \tag{HS} 
u_{txx} + 2 u_x u_{xx} + u u_{xxx} = 0 \,,
\end{equation*} 
whose solutions inevitably break down. 
We recall here a simple proof of this fact\footnote{This argument was 
worked out by K. Gleason in her Senior Thesis \cite{g}; 
see also \cite{Yin}.}. 
Consider any smooth (non-constant) periodic solution $u$ of 
the HS equation with smooth initial data $u_o$ 
(see e.g. \cite{ti}). 
Integrating (\ref{HS}) with respect to the space variable we obtain  
\begin{equation*} \label{HS1} 
u_{tx} + \frac{1}{2} u_{x}^2 + u u_{xx} = a 
\end{equation*} 
where 
$$
a = -\frac{1}{2} \int u_{x}^2 \, dx 
$$
is independent of time as is easily checked using periodicity and 
integration by parts. (Alternatively, one notices that $a$ is a multiple of the
$\dot H^1$-energy of the solution, which is conserved for HS.)
Since $u'_{o}$ is continuous we can pick a point $x_o$ where 
$u'_{o}(x_o)<0$. Then, along the characteristic $t\to x(t)$ 
defined by $x(0)=x_o$ we get an ordinary differential equation 
of the form 
$$
\frac{dw}{dt} + \frac{1}{2} w^2 = a, 
\qquad 
\mathrm{where}
\quad
w(t) = u_x(t, x(t)) 
$$
whose explicit solution 
$$
w(t) 
= 
\frac{ u'_{o}(x_o)\sqrt{-2a} - 2a \tan{(-t\sqrt{-2a}/2 )} }
{\sqrt{-2a} + u'_{o}(x_o)\tan{(t\sqrt{-2a}/2})}
$$
becomes unbounded as 
$
t \nearrow   
T_{crit} 
= 
\frac{2}{\sqrt{-2a}}\arctan{\big(\sqrt{-2a}/|u'_{o}(x_o)|\big)}. 
$
\smallskip

Note, however, that blow-up of $\mu$HS solutions can occur even if the mean 
of the initial data is not zero. More precisely, if $u$ is a smooth 
solution of (\ref{LmuHS})--(\ref{ic}) then proceeding as above 
we arrive at the equation 
$$
u_{tx} - 2\mu(u) u + \frac{1}{2} u_x^2 + u u_{xx} = r_0 \,,
$$
where the constant of integration 
$$
r_0 = -2\mu(u)^2 - \frac{1}{2} \int u_x^{2}\, dx 
$$ 
is again independent of time. 
Using Cauchy-Schwarz and conservation of the energy 
$2H_1(u) = \mu(u)^2 + \int u_x^2 dx$ 
and the mean (\ref{mean}), we then estimate 
\begin{equation*} 
\begin{split} 
u_{tx} + \frac{1}{2} u_x^2 + u u_{xx} 
&= 
r_0 + 2\mu(u)u                   \\ 
&= 
2 \mu(u_o) \Big( u(t,x) - \int u(t,x) \, dx \Big) 
- 
\frac{1}{2} \int (u'_o)^2 dx      \\ 
&\leq 
\Big(2 |\mu(u_o)| -\frac{1}{2}\|u'_o\|_{L^2} \Big)\|u'_o\|_{L^2} 
\end{split}
\end{equation*}
since the function $u - \mu(u)$ must vanish somewhere 
on the circle. The conservation law $H_1$ ensures that solutions 
remain bounded as long as they are defined. Therefore the characteristics 
are defined as long as the solution exists. 
If we arrange for the right-hand side above to be non-positive then 
along any characteristic the function $w(t)=u_x(t,x(t))$ 
will satisfy 
\begin{equation} \label{w'} 
\frac{dw}{dt} \leq -\frac{1}{2} w^2 
\end{equation} 
and so picking $x_o$ as before with $u'_o(x_o) <0$ 
we get 
$$
w(t) \leq w(0)=u'_o(x_o) <0 
$$
as long as it exists. 
Using this and solving (\ref{w'}) we therefore obtain 
for any $t \geq 0$ the estimate 
$$
w(t) \leq \frac{2u'_o(x_o)}{2 + t u'_o(x_o)} \,,
$$
from which conclude that $w(t)$ becomes infinite 
at $T_{crit}= 2/|u'_o(x_o)|$. 
We summarize these results in 
\begin{proposition} \label{wp_blowup}
For any smooth initial data satisfying 
$4|\mu(u_o)| \leq \|u'_o \|_{L^2}$ 
there exists $T_{crit}>0$ such that the corresponding solution $u(t,x)$
of the $\mu$HS equation remains bounded for $t < T_{crit}$ and satisfies 
$\|u_x(t)\|_\infty \nearrow \infty$ whenever $t \nearrow T_{crit}$. 
\end{proposition} 

Thus, the $\mu$HS equation is an example of a nonlinear PDE 
admitting solutions with finite as well as infinite 
lifespan\footnote{Results of this type for the CH equation 
can be found e.g. in \cite{ce}.}. 
As a simple illustration, the solution corresponding to 
the initial data $u_o(x) = a + \cos{2\pi x}$ blows up whenever 
$|a| \leq \pi\sqrt{2}/4$ 
while persisting indefinitely if 
$|a| \geq 4\pi^2$. 

It seems likely that Theorem \ref{wp_global} and 
Proposition \ref{wp_blowup} can be significantly strengthened 
to give the full picture of well-posedness including a complete 
characterization of the blow-up mechanism (e.g. along the lines 
of \cite{McK} or by studying the evolution of the extrema of $u_x$ 
as in \cite{ce2}). 



\section{Traveling-wave solutions} \label{waves} \nequation
In order to find traveling-wave solutions of the $\mu$HS equation 
we first regard $\mu$ as a parameter independent of $u$. 
We look for solutions  
$u(t,x) = \varphi(x-ct)$, traveling with speed $c \in \R$, 
of the equation
\begin{equation}\label{muHSparameter}
 - u_{txx}  
= 
- 2\mu u_x + 2 u_x u_{xx} + u u_{xxx},
\end{equation}
where $\mu \in \R$ is a parameter. A first analysis will 
give us a list of all traveling waves of this equation. Subsequently, by 
keeping only the solutions $\varphi$ for which the parameter $\mu$ equals 
the mean $\mu(\varphi) = \int \varphi dx$ of $\varphi$, the list will be 
narrowed down to traveling waves of the true $\mu$HS equation. 
Since $\mu(\varphi)$ is a conserved quantity, this procedure yields 
all traveling waves of $\mu$HS.

Substituting $u(t,x) = \varphi(x-ct)$ into equation (\ref{muHSparameter}) 
and integrating, we get
$$c\varphi_{xx}
   = - 2\mu\varphi+ \varphi\varphi_{xx} + \frac{1}{2}\varphi_x^2 + \frac{a}{2},
$$ 
for some constant $a \in \R$.
Rewriting this as
\begin{equation}\label{weaktrav}
           \varphi_x^2 
           = -4\mu\varphi +  \left((\varphi-c)^2\right)_{xx} + a.
\end{equation}
we have an equation that makes sense for $\varphi \in H^1_{loc}(\R)$ and 
it is natural to make the following definition: 
A function $\varphi \in H^1_{loc}(\R)$ is a {\it traveling wave
of equation (\ref{muHS})} if it satisfies (\ref{weaktrav}) in distributional 
sense for some $a \in \R$. Applying the following lemma with $v = \varphi - c$ 
we immediately deduce that any traveling-wave solution $\varphi$ is 
smooth except possibly on the boundary of the set $\{x \in \R ~|~ \varphi(x) = c\}$. A general traveling-wave solution $\varphi(x)$ is therefore smooth on the union of the (possibly countably infinite) number of open intervals that make up the complement of this boundary.
        
\begin{lemma}\label{vlemma}{\rm  (\cite{L3})}
        Let $p(t)$ be a polynomial with real coefficients. Assume that
        $v \in H^1_{loc}(\R)$ satisfies
        \begin{equation}\label{v}
           (v^2)_{xx} = v_x^2 + p(v) \quad \text{in} \quad \mathcal{D}'(\R).
        \end{equation}
         Then $v^s \in C^j(\R)$ for $k \geq 2^j$. 
       In particular, $v$ is smooth away from the boundary of 
      the set $\{x \in \R ~|~ v(x) = 0\}$.
\end{lemma}

Within each interval where $\varphi$ is smooth we may integrate (\ref{weaktrav}) once more to obtain the ODE
\begin{equation}\label{varphiODE}  
  \varphi_x^2 = \frac{-2\mu\varphi^2 +a\varphi + b}{c-\varphi},
\end{equation}  
where $b \in \R$ is another integration constant which may take on different values for different intervals.
If the roots of the polynomial $-2\mu\varphi^2 +a\varphi + b$ have nonzero imaginary part, there are no bounded solutions of (\ref{varphiODE}). Thus, we may introduce real parameter $m$ and $M$ such that $m <M$ and 
\begin{equation}\label{mM}  
  -2\mu\varphi^2 +a\varphi + b = 2\mu(M-\varphi)(\varphi - m).
\end{equation}
The next theorem follows from an analysis of (\ref{varphiODE}) along the lines of \cite{L3}. Note that (\ref{weaktrav}) is left unchanged under the substitutions
\begin{equation}\label{symmetry}
  \mu \mapsto -\mu, \quad \varphi \mapsto -\varphi, \quad c \mapsto -c,
\end{equation}  
so it is enough to consider the case $\mu > 0$.
 
\begin{theorem}\label{travth}
  Let $\mu > 0$. For a given $c \in \R$, the bounded traveling-wave solutions of speed $c$ can be classified by the values of the parameters $m, M \in \R$ as follows.
\begin{itemize}
\item[(a)] (Smooth periodic) If $m<M<c$, there is a smooth periodic
traveling wave $\varphi(x-ct)$ of (\ref{muHSparameter}) with
$\min_{x\in\R}{\varphi(x)} = m$ and $\max_{x\in\R}{\varphi(x)} =M$.

\item[(b)] (Cusped periodic) If $m<c<M$, there is a cusped periodic
traveling wave $\varphi(x-ct)$ of (\ref{muHSparameter}) with
$\min_{x\in\R}{\varphi(x)} = m$ and $\max_{x\in\R}{\varphi(x)} = c.$

\item[(c)] (Anticusped periodic) If $c<m<M$, there is an anticusped periodic
traveling wave $\varphi(x-ct)$ of (\ref{muHSparameter}) with
$\min_{x\in\R}{\varphi(x)} = c$ and $\max_{x\in\R}{\varphi(x)} = m$.

\item[(d)] (Anticusped on the line) If $c<m=M$, there is an anticusped
traveling wave $\varphi(x-ct)$ of (\ref{muHSparameter}) with
$\min_{x\in\R}{\varphi(x)} = c$ and $\lim_{x \to \pm \infty} \varphi(x) = \max_{x\in\R}{\varphi(x)} = m=M$.

\item[(e)] (Composite waves)  
Two wave-segments in $(a)-(d)$ corresponding to parameters $(m_1,M_1)$ respectively $(m_2,M_2)$ are allowed to be glued together at their crests whenever
$m_1 + M_1 = m_2 + M_2$. Any composite wave formed by joining at most countably many wave-segments in this manner is a traveling wave of (\ref{muHSparameter}).

\item[$(f)$] (Plateaus) In the special case when $m_1+M_1 = m_2+M_2 = c$, 
the composite wave may also contain plateaus, i.e. intervals on which $\varphi \equiv c$.
\end{itemize}
\end{theorem}


\subsection{Solitary waves}

From Theorem \ref{travth} we infer the following list of solitary waves (i.e. 
traveling waves with decay to zero at infinity) of 
equation (\ref{muHSparameter}). 

\begin{itemize}
\item If $\mu > 0$, then, for each $c<0$, there is an anticusped
solitary wave $\varphi(x-ct)$ with
$\min_{x\in\R}{\varphi(x)} = c$ and $\lim_{x \to \pm \infty} \varphi(x) = 0$.

\item If $\mu < 0$, then, for each $c>0$, there is a cusped
solitary wave $\varphi(x-ct)$ with
$\lim_{x \to \pm \infty} \varphi(x) = 0$ and $\max_{x\in\R}{\varphi(x)} = c $.
\end{itemize}

To verify that these are indeed anticusped respectively cusped solitary waves, 
we may determine their profiles directly. For example, choosing the case $\mu < 0$ 
so that $c > 0$, we employ (\ref{varphiODE}) with integration constants $a$ 
and $b$ set to zero to get
$$
\int \sqrt{\frac{c-\varphi}{-2\mu\varphi^2}} d\varphi = x - x_0.
$$
The integral can be computed to yield the formula
$$
\sqrt{\frac{2}{-\mu}}\biggl(\sqrt{c-\varphi} - \sqrt{c} \,\hbox{arctanh}\biggl(\sqrt{\frac{c-\varphi}{c}}\biggr)\biggr)= x - x_0.
$$
\begin{figure}[ht!]
\input epsf
\centerline{\epsfysize=0.18\vsize\epsffile{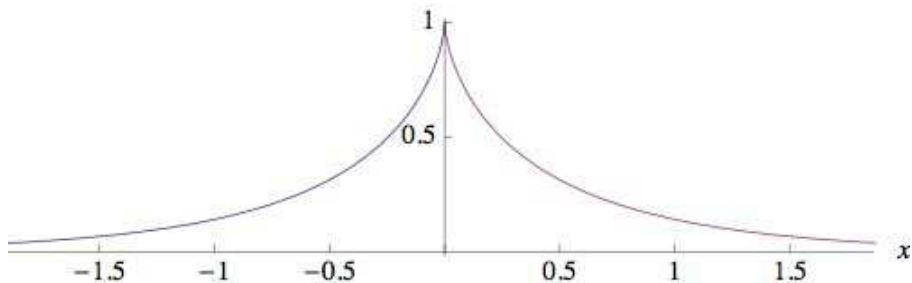}}
 \caption{Although this cusped solitary wave is a solution moving with 
speed $c = 1$ of equation (\ref{muHSparameter}) when $\mu = -1$, it is not 
a solution of the $\mu$HS equation due to the negative sign of $\mu$.}\label{solitary.pdf}
\end{figure}

Note that $x \to -\infty$ as $\varphi \searrow 0$ and that $x-x_0 \to 0$ 
as $\varphi \nearrow c$. Near the crest point we get
$$
\frac{\sqrt{2}}{3c\sqrt{-\mu}}(c-\varphi)^{3/2} + O((c-\varphi)^{5/2}) = |x- x_0|,
$$
and so, near $x = x_0$,
$$c- \varphi(x) = \biggl(\frac{3c\sqrt{-\mu}}{\sqrt{2}}\biggr)^{2/3}|x-x_0|^{2/3} 
+ O(|x-x_0|^{4/3}),$$ 
showing that there is indeed a cusp at $x = x_0$. 

Do any of these solitary waves satisfy the condition that $\mu$ 
equals $\int_\R \varphi dx$, so that they are also solutions of the 
$\mu$HS equation? The answer is no. Indeed, if $\mu > 0$, then since 
$\max_{x\in\R}{\varphi(x)} = 0$, the solitary wave has $\varphi \leq 0$ 
everywhere, and so $\mu(\varphi) = \int_\R \varphi dx \leq 0$ never equals $\mu$. 
The case $\mu < 0$ is analogous. 

\begin{theorem}
The  $\mu$HS equation does not have solitary wave solutions on the line.
\end{theorem}

There are also composite solitary wave solutions obtained by inserting cusped 
or anticusped periodic wave segments at the crests of these solitary waves 
according to (e) in Theorem \ref{travth}. One might hope that a clever 
construction of a composite wave would yields a different outcome. 
However, since $m + M = 0$ for the solitary wave and $m < M$ by definition, 
one can see that any cusped wave segments joined to the solitary wave 
must have the same sign, and thus the sign of the integral $\int_\R \varphi dx$ 
remains unchanged. 

This analysis shows that the $\mu$HS equations admits no solitary waves, 
and hence no solitons, on the line. In the periodic setting for $\mu$HS, 
however, we will see that traveling waves which are natural candidates for solitons do exist. 

\bigskip


\subsection{Traveling waves on the circle}
In view of the symmetry (\ref{symmetry}) we restrict ourselves to 
the case $\mu > 0$ in this section.

Consider equation (\ref{varphiODE}) 
\begin{equation}\label{varphiODEmM}  
  \varphi_x^2 = \frac{2\mu(M-\varphi)(\varphi - m)}{c-\varphi},
\end{equation}
written in terms of $m$ and $M$ according to (\ref{mM}). By (a), (b), and (c) 
of Theorem \ref{travth} there is a unique periodic traveling wave 
$\varphi_{m,M,\mu}$ of equation (\ref{muHSparameter}) associated to 
each point $(m,M,\mu)$ such that (a) $m<M<c$, (b) $m < c < M$, or (c) 
$c < m < M$. By enforcing the constraint 
\begin{equation}\label{constraint}  
  mean(\varphi_{m,M,\mu}) = \mu
\end{equation}
we obtain the subset of these that are also solutions of the $\mu$HS equation.
From the way $\mu$ appears in (\ref{varphiODEmM}) we infer that (see (\ref{meanvarphimMmu}) below)
\begin{equation}\label{muscaling}  
  mean(\varphi_{m,M,\mu}) = \frac{1}{\sqrt{\mu}} mean(\varphi_{m,M,1}).
\end{equation}
The constraint (\ref{constraint}) therefore translates to 
$mean(\varphi_{m,M,1}) = \mu^{3/2}$
and is seen to admit a unique solution $\mu > 0$ exactly when 
$mean(\varphi_{m,M,1}) > 0$. Whereas this condition fails 
to hold in the case of solitary waves as noted above, it is satisfied 
for large classes of periodic waves. Rather than finding the most general 
choices of $m$, $M$, and $c$ for which it is fulfilled, we simply note 
that if the minimum of $\varphi_{m,M,1}$ is positive then its mean is 
certainly also positive. Applying Theorem \ref{travth} with $\mu = 1$ 
we can therefore deduce existence of the following classes of traveling 
waves parametrized by $m, M, c \in \R$.

\begin{itemize}
\item[(a)] (Smooth periodic) Whenever $0 < m<M<c$, there is a smooth periodic
traveling wave $\varphi(x-ct)$ of $\mu$HS with
$\min_{x\in\R}{\varphi(x)} = m$ and $\max_{x\in\R}{\varphi(x)} =M$.

\item[(b)] (Cusped periodic) Whenever $0 < m<c<M$, there is a cusped periodic
traveling wave $\varphi(x-ct)$ of $\mu$HS with
$\min_{x\in\R}{\varphi(x)} = m$ and $\max_{x\in\R}{\varphi(x)} =c$.

\item[(c)] (Anticusped periodic) Whenever $0 < c<m<M$, there is 
an anticusped periodic
traveling wave $\varphi(x-ct)$ of $\mu$HS with
$\min_{x\in\R}{\varphi(x)} = c$ and $\max_{x\in\R}{\varphi(x)} =m$.
\end{itemize}

In the rest of this section we will derive explicit expressions 
(involving elliptic integrals) for the period and the mean of these waves. 
As a consequence of the analysis below, we will see that
there are plenty of traveling waves of period one, i.e. on the circle $S^1=\R/\Z$.

\begin{theorem}\label{periodoneTh}
For any $c \neq 0$ the $\mu$HS equation 
has a one-parameter family of smooth periodic traveling waves of period one and
a one-parameter family of cusped periodic traveling waves of period one
moving with speed $c$.
\end{theorem}

\subsubsection{Smooth periodic waves}
From (\ref{varphiODEmM}) it follows that
\begin{equation}\label{xminusx0}
  x - x_0 = \frac{1}{\sqrt{ 2 \mu}}\int_{\varphi_0}^\varphi \sqrt{\frac{c-\varphi}{(M-\varphi)(\varphi - m)}}d\varphi.
\end{equation}  
Substituting $\varphi = m + (M-m)\sin^2\theta$, we find
$$x - x_0 = \sqrt{\frac{2(c-m)}{\mu}}\int_{\theta_0}^\theta \sqrt{1 - \frac{M-m}{c-m}\sin^2\theta} d\theta= \sqrt{\frac{2(c-m)}{\mu}} E\left(\theta
   \left|\frac{M-m}{c-m}\right.\right)\biggr|_{\theta_0}^\theta,$$
which expresses $x-x_0$ in terms of the incomplete elliptic integral of the second kind
$$E(\phi | r) := \int_0^\phi \sqrt{1- r\sin^2\theta}d\theta, \qquad -\frac{\pi}{2} < \phi < \frac{\pi}{2}.$$

Now for smooth periodic waves the period is given by twice the integral from the minimum $\varphi_0 = m$ to the maximum $\varphi_1 = M$, or, in terms of $\theta$, from $\theta_0 = 0$ to $\theta_1 = \frac{\pi}{2}$. We deduce that 
\begin{equation}\label{periodmMmu}
  period(\varphi_{m,M,\mu}) = 2\sqrt{\frac{2(c-m)}{\mu}} E\left(\frac{M-m}{c-m}\right)
\end{equation}  
where $E(r) := E(\pi/2 | r)$ is the complete elliptic integral of the second kind.

Similarly, we compute
$$
\int_{x_0}^x \varphi dx = \frac{1}{\sqrt{ 2 \mu}}\int_{\varphi_0}^\varphi \varphi \sqrt{\frac{c-\varphi}{(M-\varphi)(\varphi - m)}}d\varphi
$$
$$
= \sqrt{\frac{2(c-m)}{\mu}}\int_{\theta_0}^\theta (m + (M-m)\sin^2\theta)\sqrt{1 - \frac{M-m}{c-m}\sin^2\theta} d\theta\,.
$$
Then integration over a full period yields
\begin{equation}\label{meanvarphimMmu}  
  mean(\varphi_{m,M,\mu}) = \frac{2}{3}\sqrt{\frac{2(c-m)}{\mu}} \left((2 (m+M) - c) E\left(\frac{M-m}{c-m}\right)+(c - M)
   K\left(\frac{M-m}{c-m}\right)\right),
\end{equation}   
where 
$$K(r) := \int_0^{\pi/2} (1 - r \sin^2 \theta)^{-1/2} d\theta,$$
is the complete elliptic integral of the first kind.
Observe that the mean indeed scales with $\mu$ in accordance with (\ref{muscaling}). 

Solving the constraint (\ref{constraint}) for $\mu$ we see that the mean of the traveling wave $\varphi$ of the $\mu$HS equation associated to the point $(m,M)$ is
$$\mu(\varphi) = \left(\frac{2}{3}\sqrt{2(c-m)} \left((2 (m+M) - c) E\left(\frac{M-m}{c-m}\right)+(c - M)K\left(\frac{M-m}{c-m}\right)\right)\right)^{2/3}$$ 
and substituting this value of $\mu$ into (\ref{periodmMmu}), we find that the corresponding period is
\begin{equation}\label{periodvarphi}  
  period(\varphi) = \frac{2\cdot 3^{1/3}\sqrt{c-m}E\left(\frac{M-m}{c-m}\right)}{\left(\sqrt{c-m} \left((2 (m+M) - c) E\left(\frac{M-m}{c-m}\right)+(c - M)K\left(\frac{M-m}{c-m}\right)\right)\right)^{1/3}}.
\end{equation}

\begin{figure}[ht!]
\input epsf
\centerline{\epsfysize=0.33\vsize\epsffile{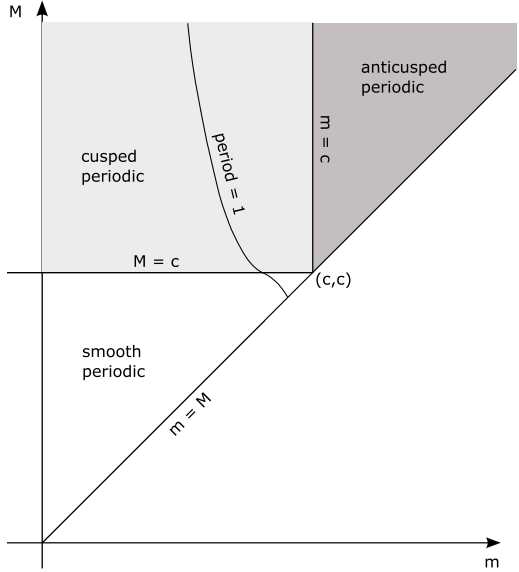}}
 \caption{ For each $c \neq 0$ there are one-parameter families of both smooth and cusped period-one traveling waves of the $\mu$HS equation.}\label{muHStravelingwaves.pdf}
\end{figure}

\subsubsection{Cusped periodic waves}
For cusped periodic waves ($m < c < M$) the period is given by twice the integral (\ref{xminusx0}) from $\varphi_0 = m$ to $\varphi_1 = c$, i.e. from $\theta_0 = 0$ to $\theta_1 = \arcsin\left(\sqrt{\frac{c-m}{M-m}}\right)$. We get
$$period(\varphi_{m,M,\mu}) = 2\sqrt{\frac{2(c-m)}{\mu}} E\left(\arcsin\left(\sqrt{\frac{c-m}{M-m}}\right)
   \left|\frac{M-m}{c-m}\right.\right).$$
Similarly, formula (\ref{meanvarphimMmu}) for $mean(\varphi_{m,M,\mu})$ and formula (\ref{periodvarphi}) for the period also hold in the cusped case with the replacements
$$E\left(\frac{M-m}{c-m}\right) \to E\left(\arcsin\left(\sqrt{\frac{c - m}{M - m}}\right)\biggl|\frac{M-m}{c-m}\right),$$
and
$$
K\left(\frac{M-m}{c-m}\right) \to F\left(\arcsin\left(\sqrt{\frac{c-m}{M-m}}\right)\biggl|\frac{M-m}{c-m}\right).
$$
where 
$$
F(\phi | r) := \int_0^\phi (1 - r \sin^2 \theta)^{-1/2} d\theta,
$$
is the incomplete analog of $K$. An analysis of the expression for the period 
in (\ref{periodvarphi}), together with its cusped analog, shows that 
the level set of period-one solutions looks qualitatively as in 
Figure \ref{muHStravelingwaves.pdf}. In particular, we see that 
for each $c > 0$ there are one-parameter families of both smooth and cusped 
period-one traveling-wave solutions of the $\mu$HS equation. 
This establishes Theorem \ref{periodoneTh}.

\bigskip


\section{Geometry of the group of circle diffeomorphisms} 
\label{curvature}\nequation
The objective of this section is to compute the curvature of 
$\mathcal{D}^s(S^1)$ endowed with the $\mu$HS metric (\ref{Metric}). 
We establish a local formula for a connection compatible with the metric
and use it to obtain an exact expression for the curvature tensor.
In particular, we prove that sectional curvatures in both vertical and horizontal
two-dimensional directions on  $\mathcal{D}^s(S^1)$ are positive.


\subsection{Christoffel symbols and the covariant derivative}
Define a symmetric bilinear map $\Gamma_{id}( \cdot, \cdot)$ 
on $T_{id} \mathcal{D}^s(S^1) \simeq H^s(S^1)$ by 
\begin{equation}\label{Gammadef}
  \Gamma_{id}( u, v) 
= - A^{-1}\partial_x\bigl(\mu(u)v 
+ \mu(v)u + \frac{1}{2}u_x v_x\bigr), \qquad u,v \in H^s(S^1),
  \end{equation}  
and extend it to a bilinear map 
$\Gamma_\eta( \cdot, \cdot):T_\eta \mathcal{D}^s(S^1) 
\times T_\eta \mathcal{D}^s(S^1)  \to T_\eta \mathcal{D}^s(S^1)$ 
for any $\eta \in \mathcal{D}^s(S^1)$ 
using right invariance. 

Next, define the operator $\nabla$ locally in a chart at $\eta$ 
by the formula 
\begin{equation}\label{covderlocal}
  (\nabla_X Y)(\eta) = DY(\eta) \cdot X(\eta) 
- 
\Gamma_\eta( Y(\eta), X(\eta)) 
\end{equation}
for any vector fields $X$ and $Y$ on $\mathcal{D}^s(S^1)$. 
It turns out that $\nabla$ is the smooth Levi-Civita connection 
of the right-invariant $\mu$HS metric whose Christoffel symbols 
are given by the map $\Gamma$. 
The proof of the following proposition is straightforward 
(see for example \cite{em}).

\begin{proposition}\label{covderprop}
The bilinear map $\nabla$ defined by (\ref{covderlocal}) is 
the unique torsion-free covariant derivative on $\mathcal{D}^s(S^1)$ 
compatible with the right-invariant metric (\ref{Metric}).
\end{proposition}
\medskip


\subsection{Sectional curvature}
\label{curvsec}
In order to compute the sectional curvature of $\mathcal{D}^s(S^1)$ 
equipped with the $\mu$HS metric $\langle \cdot, \cdot \rangle$ 
defined in (\ref{Metric}) recall first that 
its curvature tensor $R$ can be determined from the formula 
\begin{equation}\label{RXYZ}
  R(X,Y)Z = \nabla_X \nabla_Y Z - \nabla_Y \nabla_X Z - \nabla_{[X, Y]} Z 
\end{equation}  
for any vector fields $X,Y$ and $Z$ on $\mathcal{D}^s(S^1)$. 
By right invariance it will be sufficient to compute $R$ at the identity. 
We will write $\Gamma(u,v)$ to denote the Christoffel symbol 
$\Gamma_{id}(u, v)$.

\begin{theorem}\label{curvth}
For any $u,v \in T_{id} \mathcal{D}^s(S^1) \simeq H^s(S^1)$, it holds 
that
\begin{align}\label{curv}
 \langle R(u, v)v, u \rangle 
= & \langle \Gamma(u,v), \Gamma(u,v)\rangle
- \langle \Gamma(u,u), \Gamma(v,v)\rangle -3\mu( u_x v )^2
		\nonumber \\
= & \mu(u)^2\bigl(\mu(v^2) + \mu(v_x^2)\bigr) 
           + \mu(v)^2\bigl(\mu(u^2) + \mu(u_x^2)\bigr) 
		\\ \nonumber
&+ \mu(u)\mu( (vu_x - uv_x)v_x)
+ \mu(v)\mu((uv_x - vu_x)u_x )
		\\ \nonumber
&-2\mu(u)\mu(v)\bigl(\mu(uv) + \mu(u_xv_x)\bigr)
		\\ \nonumber
&- \frac{1}{4}\left(\int u_x v_x \,dx\right)^2 
  + \frac{1}{4} \int u_x^2\,dx \int  v_x^2\,dx 
   - 3\mu( u_x v )^2\,.
\end{align}
\end{theorem}

\proofbegin ~ 
For any $u, v$ and $w \in T_{id}^s\mathcal{D}(S^1) \simeq H^s(S^1)$ 
the local formula corresponding to the expression in (\ref{RXYZ}) 
reads 
\begin{align*}
  R(u,v)w = &D_1\Gamma_{id}(w,u)v - D_1\Gamma_{id}(w,v)u
  		 \\
 &  + \Gamma_\eta(\Gamma_{id}( w,v) ,u) - \Gamma_{id}(\Gamma_\eta  (w,u),v)
\end{align*}
where $D_1$ denotes differentiation of $\Gamma_\eta$ with respect to $\eta$:
$D_1\Gamma_{id}(w,u)v = \frac{d}{d\epsilon}|_{\epsilon = 0} \Gamma_{id + \epsilon v}(w,u)$.
We compute
$$
D_1\Gamma_{id}( w, u)v
= -\Gamma_{id}( w_xv, u) - \Gamma_{id}( w, u_xv) + \Gamma_{id}( w,u)_x v,$$
and so
\begin{align*}
\langle R(u, v)v, u \rangle 
 =&  \bigl\langle -\Gamma(v_xv, u) - \Gamma(v, u_xv) + \Gamma(v_xu, v) + \Gamma(v, v_xu), u \bigr\rangle
		\\
& + \bigl\langle \Gamma(v,u)_x v - \Gamma(v,v)_x u + \Gamma(\Gamma(v,v),u) - \Gamma(\Gamma(v,u),v), u \bigr\rangle.
\end{align*}
Observe that
$$
\Gamma(u, v) 
= 
\frac{1}{2}\bigl((uv)_x - ad^*_v(u) - ad^*_u(v)\bigr),
$$
where $ad^*_v(u) = A^{-1}(2v_xAu + vAu_x)$, see (\ref{Bdef}).
Furthermore, note that the last four terms in the expression 
for $\langle R(u, v)v, u \rangle$ 
can be rewritten as 
$$
\langle \Gamma(v,u)_x  - \Gamma(v,v)_x u + \Gamma(\Gamma(v,v),u) - \Gamma(\Gamma(v,u),v), u \rangle
$$
$$
= \langle \Gamma(u,v), \Gamma(u,v)\rangle
- \langle \Gamma(u,u), \Gamma(v,v)\rangle
$$
$$ 
- \langle u_xv, \Gamma(v, u)\rangle
+ \langle uu_x, \Gamma(v,v)\rangle.
$$
Therefore
\begin{equation*}
\begin{split} 
\langle R(u, v)v, u \rangle 
= &\langle -\Gamma(v_xv, u) - \Gamma(v, u_xv) + \Gamma(v_xu, v) + \Gamma(v, v_xu), u \rangle 
	\\
&+ \langle \Gamma(u,v), \Gamma(u,v)\rangle - \langle \Gamma(u,u), \Gamma(v,v)\rangle
	\\
& - \langle u_xv, \Gamma(v, u)\rangle
+ \langle uu_x, \Gamma(v,v)\rangle.
\end{split}
\end{equation*}

Now, using the formula for the Christoffel symbols in \eqref{Gammadef} 
a straightforward computation gives 
\begin{equation*} 
\begin{split} 
\bigl\langle -\Gamma(v_xv, u) 
- 
\Gamma(v, u_xv) 
+ 
\Gamma(v_xu, v) + \Gamma(v, v_xu), u \bigr\rangle   
&- 
\langle u_xv, \Gamma(v, u)\rangle
+ 
\langle uu_x, \Gamma(v,v)\rangle   \\ 
&= 
-3\biggl(\int u_xv dx\biggr)^2.
\end{split}
\end{equation*} 
This establishes the first curvature formula of Theorem \ref{curvth}. 

To obtain the second equality in \eqref{curv} we simply substitute 
the expression for the Christoffel symbol $\Gamma$ into
$$ 
\langle \Gamma(u,v), \Gamma(u,v)\rangle
- \langle \Gamma(u,u), \Gamma(v,v)\rangle -3\biggl(\int u_xv dx\biggr)^2
$$
and integrate by parts to get
$$
-\int (\mu(u)v + \mu(v)u + \frac{1}{2}u_xv_x) \partial_x A^{-1}\partial_x(\mu(u)v + \mu(v)u + \frac{1}{2}u_xv_x)dx
$$
$$
+ \int (2\mu(u)u + \frac{1}{2}u_x^2) \partial_xA^{-1}\partial_x(2\mu(v)v 
+ \frac{1}{2}v_x^2)dx -3\biggl(\int u_xv dx\biggr)^2.
$$
Using the identity 
$\partial_xA^{-1}\partial_xf = A^{-1}\partial_x^2 = -f + \mu(f)$  
we simplify this expression and obtain the right-hand side 
of (\ref{curv}).
\proofend

From Theorem \ref{curvth} we can immediately deduce an expression for 
the sectional curvature $K(u,v)$ of the plane spanned by two vectors 
$u,v \in T_{id} \mathcal{D}^s(S^1)$. 
For this purpose we may assume, after taking linear combinations, 
that $u$ and $v$ are orthonormal with respect to 
the $\mu$HS metric 
$\langle \cdot, \cdot \rangle$ and that $\mu(v) = 0$, 
i.e.
$$
\mu(u)^2 + \int_{S^1} u_x^2 \, dx = 1, 
\quad 
\mu(v) = 0, 
\quad 
\int_{S^1} v_x^2 \, dx = 1, 
\quad 
\text{and}
\quad 
\int_{S^1} u_x v_x \,dx = 0.
$$
With these assumptions, Theorem \ref{curvth} yields
\begin{equation}\label{Suv}
K(u,v) 
=  
\mu(u)^2\bigl(\mu(v^2) + 1\bigr) 
+ 
\mu(u)\mu\big( (vu_x - uv_x)v_x \big) 
+ 
\frac{1}{4} \int u_x^2 \,dx  - 3\mu( u_x v )^2\,.
\end{equation}

\begin{proposition}\label{prop:paralH}
For any orthonormal vectors $u,v \in T_{id} \mathcal{D}^s(S^1)$ with $\mu(u) = \mu(v) = 0$ the sectional curvature $K(u,v)$ of the plane spanned by $u$ and $v$ satisfies
\begin{equation}\label{Sinequality}  
  K(u,v) = \frac{1}{4} - 3\mu( u_x v )^2 \geq \frac{1}{4}\left(1 - \frac{3}{\pi^2}\right) >0.
\end{equation}
\end{proposition}
\proofbegin ~ 
In this case (\ref{Suv}) reduces to
\begin{equation}\label{S14}  
  K(u,v) =  \frac{1}{4} - 3\mu( u_x v )^2\,.
\end{equation}
By the Cauchy-Schwartz inequality it holds that
$$\left(\int u_x v dx\right)^2 \leq \int u_x^2 dx \int v^2 dx = \int v^2 dx.$$
Since $\mu(v) = 0$ we also have 
$$\int v^2 dx \leq \frac{1}{4\pi^2} \int v_x^2 dx = \frac{1}{4\pi^2},$$
so that (\ref{Sinequality}) follows from (\ref{S14}).
\proofend

\begin{remark}
{\rm
The origin of the term $- 3\mu( u_x v )^2$ in (\ref{Sinequality}) is as follows. For any $u,v\in T_{id} \mathcal{D}^s(S^1)$ decomposed as $u = \tilde{u} + \mu(u)$ and $v = \tilde{v} + \mu(v)$ with $\mu(\tilde{u}) = \mu(\tilde{v}) = 0$, Theorem \ref{curvth} yields 
\begin{align}\label{curvaturedecomposed}
\langle R(u,v)v, u \rangle
& =
\text{$\mu$-terms}
+
\langle R^{HS}(\tilde{u}, \tilde{v})\tilde{v}, \tilde{u} \rangle
- 3\mu( u_x v )^2,
\end{align}
where `$\mu$-terms' stands for those terms containing either 
$\mu(u)$ or $\mu(v)$ (or both) as a factor, and 
$$
\langle R^{HS}(\tilde{u}, \tilde{v})\tilde{v}, \tilde{u} \rangle
= -
\frac{1}{4} \left(\int_{S^1} \tilde{u}_x \tilde{v}_x dx \right)^2
+\frac{1}{4} \int_{S^1} \tilde{u}_x^2 dx \int_{S^1} \tilde{v}_x^2 dx
$$
is the curvature of the $\dot{H}^1$ metric (\ref{metric}) 
corresponding to the HS equation, see \cite{L}. 
Then 
the middle term in (\ref{curvaturedecomposed}) corresponds to the curvature of
the base $\mathcal{D}^s(S^1)/S^1$ contributing the constant term $\frac{1}{4}$ in (\ref{Sinequality}), while the last term can be seen
as the correction term in the general O'Neill formula \cite{on} 
$$
\langle R^M(X, Y)Y, X \rangle
=
\langle R^B(\tilde{X}, \tilde{Y})\tilde{Y}, \tilde{X} \rangle
- \frac{3}{4} \| [\tilde{X}, \tilde{Y}]^V \|^2 \,,
$$
for a Riemannian submersion $M \to B$, where 
$\tilde{X}, \tilde{Y}$ are vectors tangent to the base $B$, 
and $X, Y$ are their horizontal lifts to $M$.
Indeed, for horizontal vectors $u$ and $v$ the $\mu$-terms in 
(\ref{curvaturedecomposed}) vanish
so that since $\mu(u_x v)^2 = \mu([\tilde{u},\tilde{v}])^2/4$ the
identification is clear. 
}
\end{remark}

To summarize, let us  decompose the tangent space at the identity
$T_{id} \mathcal{D}^s(S^1) = H \oplus V$ in such a way that 
$H$ consists of all mean zero functions on $S^1$ and $V \simeq \R$ 
are the constants
(so that $u = \tilde{u} + \mu(u)$ with $\mu(\tilde{u})=0$).
Proposition \ref{prop:paralH} shows that the sectional curvature 
for any plane contained in (i.e. parallel to) the subspace $H$ is 
strictly positive. Now turn to the planes perpendicular to $H$, 
i.e. planes containing constant functions, which constitute $V$. 

\begin{proposition}\label{S1vprop}
Let $v \in T_{id} \mathcal{D}^s(S^1)$ be orthonormal to the constant function $1$. Then the sectional curvature of the plane spanned by $1$ and $v$ is
$$K(1,v) = \int v^2 dx >0.$$
\end{proposition}
\proofbegin ~ 
Apply formula (\ref{Suv}) with $u \equiv 1$ and $v$ satisfying $\mu(u) = 1$, $u_x \equiv 0$, $\mu(v) = 0$, and $\int_{S^1} v_x^2 dx = 1$.
\proofend

\begin{remark}
{\rm Applying Proposition \ref{S1vprop} to the functions $\{v_k(x) = \frac{\sqrt{2}\sin{kx}}{k}\}_{k=1}^\infty$, we obtain a sequence of two-dimensional subspaces whose sectional curvatures $K(1, v_k) = \frac{1}{k^2}$ approach zero as $k \to \infty$.}
\end{remark}

\begin{remark}
{\rm
The group $\mathcal{D}^s(S^1)$ can be thought of as a solid torus, with a
cross-section isomorphic to $\mathcal{D}^s(S^1)/S^1$. As we mentioned 
in the introduction, there is a natural 
Riemannian submersion of the metric  $\langle \cdot, \cdot \rangle$ on 
$\mathcal{D}^s(S^1)$ to the $\dot H^1$-metric on $\mathcal{D}^s(S^1)/S^1$ 
defined in (\ref{metric}).
The latter is isometric to a (contractible) piece $U\subset S^\infty$
of the unit sphere $S^\infty$ in $L^2(S^1)$, see \cite{L}.\footnote{One 
can show, however, that there is no isometric embedding of the 
group 
$\mathcal{D}^s(S^1)$ endowed with the $\mu$HS metric 
to $S^\infty\times S^1$, equipped with the
direct product metric.} This 
set $U$ has a finite diameter, and any geodesic on it reaches the boundary
of $U$ before acquiring a conjugate point.

This geometric picture sheds light on the (non)existence of global solutions 
proved in the section  on well-posedness. 
Since the group $\mathcal{D}^s(S^1)$ 
with the metric  $\langle \cdot, \cdot \rangle$ is a solid torus 
whose cross-section is a slight deformation of the open
subset $U$ (of finite diameter), it is natural to expect 
that geodesics starting almost parallel to
the cross-section (i.e. with small values $\mu(u)$) would run into the
boundary of the solid torus (and only deviate slightly from the HS
geodesics, which hit the boundary in short time). On the other hand, 
geodesics starting `perpendicular' to the base would wind around the torus
indefinitely and correspond to global solutions. This is exactly
the picture which was made precise in Section \ref{WP}.
}
\end{remark}


\section{$\mu$HS as a bi-hamiltonian system on the Virasoro group}
\label{Virasoro}\nequation

In Section \ref{GeoSec} we showed that the $\mu$HS equation 
is the equation of 
the geodesic flow of the right-invariant metric 
(\ref{Metric}) on the diffeomorphism group $\mathcal{D}^s(S^1)$, 
see Theorem \ref{Eulerth}. In other words, the $\mu$HS equation is
the Euler equation on the Lie algebra\footnote{Strictly speaking, 
$T_\mathrm{id}\mathcal{D}^s(S^1)$ may be considered a Lie algebra only 
when $s=\infty$. The reader can either make this assumption throughout 
the section or else take the Sobolev index $s$ sufficiently large 
for all constructions to be rigorously justified.} 
$\mathfrak{vect} = T_\mathrm{id}\mathcal{D}^s(S^1)$. 
By using  the inner product (\ref{metric}) to identify 
$\mathfrak{vect}$ with its dual algebra $\mathfrak{vect}^\ast$, 
the equation can be viewed as Hamiltonian 
with respect to the canonical Lie-Poisson structure 
and the Hamiltonian $H_1$. 

In this section we shall carry out this construction on the 
Virasoro group, which is the universal 
central extension of $\mathcal{D}(S^1)$ or, equivalently, on 
the corresponding Lie algebra called the Virasoro algebra. 
As with the KdV and CH equations the Virasoro group and algebra 
provide a natural setting 
for a geometric interpretation of the bi-Hamiltonian structure 
of the $\mu$HS equation discussed in Section \ref{Biham}. 
We will see that passing to the Virasoro algebra  introduces an extra term 
$2ku_x$. The corresponding Euler equation 
reads\begin{equation}\label{muHSwithux}
- u_{txx} = - 2\mu(u) u_x + 2ku_x + 2 u_x u_{xx} + u u_{xxx} 
\qquad 
(k \in \R). 
\end{equation} 
The situation is analogous to that of the CH equation
which on the Virasoro algebra becomes (cf. \cite{K-M})
$$  
u_t-u_{txx}+2k u_x + 3uu_x = 2u_xu_{xx}+uu_{xxx} 
\qquad 
(k \in \R).
$$
In this case the parameter $k$ is related to the critical shallow 
water speed, see \cite{J1}.\footnote{One should also mention 
that in the case of the CH equation there is an important 
difference in regularity properties of the associated Riemannian 
exponential maps on the Virasoro and the diffeomorphism groups, 
see \cite{CKKT}.}

Recall that the Virasoro algebra is the vector space 
$\mathfrak{vir} = \mathfrak{vect} \oplus \R$ 
endowed with the following commutator between pairs 
$$
[(v,b),(w,c)]=\left(v_x w - v w_x, \int v_x w_{xx} \, dx \right)
$$
where $v,w \in \mathfrak{vect}$ and $b,c \in \R$. 
Define an inner product on this algebra by 
\begin{equation} \label{metricVir}
\langle  (v,b) ,(w, c)\rangle
=
\langle  v,w \rangle + bc 
\end{equation}
where $\langle \cdot,\cdot \rangle$ is the inner product 
in (\ref{metric}), 
and use it to identify the dual $\mathfrak{vir}^*$ with 
$\mathfrak{vir}$ by the pairing 
$$
\langle (m,a), (v, b) \rangle^* = \int m v \, dx + ab 
$$ 
where $(u,a) \mapsto (Au,a)=(m,a)$ is the corresponding inertia operator. 

Let $h, f: \mathfrak{vir}^\ast \to \R$ be arbitrary smooth functions. 
The dual space $\mathfrak{vir}^\ast$ carries the canonical 
{\it Lie-Poisson bracket} 
\begin{equation} \label{LP} 
\{h, f\}(m,a) 
= 
\Big\langle (m, a), \Big[
\Big(\frac{\delta h}{\delta m},\frac{\delta h}{\delta a}\Big), 
\Big(\frac{\delta f}{\delta m},\frac{\delta f}{\delta a}\Big) 
\Big]\Big\rangle^*.
\end{equation} 
On the other hand, given any point $(m_0,a_0)$ in $\mathfrak{vir}^\ast$, 
we can also define 
\begin{equation} \label{frozen} 
\{h, f\}_0(m,a) 
= 
\Big\langle (m_0, a_0), \Big[
\Big(\frac{\delta h}{\delta m},\frac{\delta h}{\delta a}\Big), 
\Big(\frac{\delta f}{\delta m},\frac{\delta f}{\delta a}\Big)
\Big]\Big\rangle^* 
\end{equation} 
called a {\em frozen} (or constant) Poisson bracket. 
It is not difficult to check that both brackets form a Poisson pair 
in the sense that a linear combination of these brackets is again 
a Poisson bracket. 

Our main result in this section is 
\begin{theorem} \label{thm:biH} 
For any $k \in \R$ the $\mu$HS equation (\ref{muHSwithux}) is Hamiltonian 
with respect to the two Poisson structures (\ref{LP}) and (\ref{frozen}) 
on $\mathfrak{vir}^*$ where the second bracket is frozen at the point 
$(m_0,a_0)=(0,1)$. 
\end{theorem}

\proofbegin~ 
From the relations
\begin{align*} 
\big\langle ad^*_{(v,b)}(m,a), (w, c) \big\rangle^*
&= \big\langle (m,a), [(v,b), (w, c)] \big\rangle^*
	\\
&= \Big\langle (m,a), 
\Big(v_xw - v w_x, \int v_x w_{xx} \, dx \Big) \Big\rangle^*
\end{align*}
holding for any elements $(m,a) \in \mathfrak{vir}^*$ 
and 
$(v,b), (w,c) \in \mathfrak{vir}$, 
we deduce an explicit form of the coadjoint action 
on $\mathfrak{vir}^\ast$ 
\begin{equation}\label{adstar}  
  ad^*_{(v,b)}(m,a) = (m_xv + 2mv_x + av_{xxx}, 0).
\end{equation}
Since the Euler equation on $\mathfrak{vir}^\ast$ is of the form 
$(m,a)_t = -ad^*_{(u,a)}(m,a)$ 
we obtain the system 
$$
\begin{cases}
	& m_{t} = -m_xu - 2mu_x - au_{xxx} \\
	& a_t = 0. 
\end{cases}
$$
The second of the above equations immediately implies that 
$a$ must be constant and using $m = Au$ the first equation 
can be rewritten as 
$$
-u_{txx} = - 2\mu(u) u_x + 2 u_x u_{xx} + u u_{xxx} - au_{xxx}. 
$$
Substituting $u \to u + a$ we recover the $\mu$HS equation in 
the more familiar form (\ref{muHSwithux}) with the constant $k = -a$. 
Since the Euler equation on $\mathfrak{vir}^\ast$ is automatically 
Hamiltonian with respect to the Lie-Poisson structure (\ref{LP}) 
and with 
$$
H_1(m,a) = \frac{1}{2}\Big( \mu(u)^2 + \int u_x^2 \, dx + a^2 \Big) 
$$ 
as its Hamiltonian function, the first part of the theorem 
follows. 

To complete the proof, observe that the corresponding Hamiltonian equation 
for the frozen Poisson bracket (\ref{frozen}) and any function 
$h: \mathfrak{vir}^\ast \to \R$ reads 
$$
(m,a)_t
= 
-ad^\ast_{\left(\frac{\delta h}{\delta m},\frac{\delta h}{\delta a}\right)} 
(m_0, a_0), 
$$
which, with the help of the formula (\ref{adstar}) for 
the coadjoint action, gives 
\begin{equation} \label{system} 
\begin{cases}
	& m_{t} = -m_{0x} \dfrac{\delta h}{\delta m} 
                  - 2m_0 \Big(\dfrac{\delta h}{\delta m}\Big)_x  
                  - a_0 \Big(\dfrac{\delta h}{\delta m}\Big)_{xxx}   \\
	& a_t = 0. 
\end{cases}
\end{equation} 
Next, choose the Hamiltonian function $h$ to be 
$$
H_2(m,a)= \int \Big(\mu(u)u^2 +\frac{1}{2}uu_x^2 + kmu\Big) dx 
$$
and compute its variational derivative
$$
\frac{\delta H_2}{\delta m} 
= 
A^{-1}\left(\mu(u^2)+2\mu(u)u -\frac{1}{2}u_x^2 -uu_{xx}+km\right),
$$
cf. equation (\ref{variationalderivatives}).
Substituting this expression into (\ref{system}) and picking 
$(m_0, a_0) = (0, -1)$ gives 
$$
\begin{cases}
	& m_{t} = \partial_x^3A^{-1}\left(\mu(u^2) + 2\mu(u)u 
                     - \frac{1}{2}u_x^2 - u u_{xx} + km\right), \\
	& a_t = 0.
\end{cases}
$$ 
Since $\partial_x^3A^{-1} = -\partial_x$ and 
$m_t = -u_{txx}$, 
the first of the above equations becomes 
$$
-u_{txx}= - 2\mu(u)u_x + 2u_xu_{xx} + u u_{xxx} + ku_{xxx} 
$$
and a substitution $u \to u - k$ gives equation (\ref{muHSwithux}).
\proofend

\begin{remark}
{\rm 
According to the classification given in \cite{K-M}, on 
the Virasoro dual there are 
generically (i.e. at each point of the `highest-dimensional' 
Virasoro coadjoint  orbits) only three types of Poisson pairs 
consisting of the Lie-Poisson and frozen brackets. 
Each such pair is represented by one of the three integrable 
equations: KdV, CH and HS. 
On the other hand, Theorem \ref{thm:biH} shows that the $\mu$HS equation 
has the same bihamiltonian structure as KdV and corresponds to 
the same Poisson pair. 
How does this new integrable equation fit in the classification 
of \cite{K-M}? 


Observe that a Poisson pair determines a bihamiltonian equation 
only up to a choice of a Casimir (see Corollary 6.4 of \cite{K-M}) 
and this choice is different for the $\mu$HS and the KdV equations. 
In fact, it is possible to choose Casimirs at each step of the 
%
%
Lenard-Magri scheme. 
In the case of $\mu$HS and KdV we start with the same Casimir 
$\int u dx$ (cf. $H_0$ in Section \ref{manyham})
but the subsequent choices of Casimirs 
are different, which  results in two different
integrable equations albeit having the same bihamiltonian structure.
}
\end{remark}

\medskip
\begin{ackn}
B.K. is grateful to the IHES in Bures-sur-Yvette and the MPI in Bonn
for kind hospitality.  
This research was partially supported by an NSERC research grant. 
\end{ackn}

\bigskip



\end{document}